\documentclass[reqno, 12pt]{amsart}
\usepackage{epsfig,amsmath,amsfonts,latexsym}
\usepackage{amsthm}
\usepackage[abs]{overpic}
\usepackage[usenames,dvipsnames]{color}
\usepackage{palatino}
\usepackage[hyperfootnotes=false]{hyperref}
\usepackage{caption}
\usepackage{dsfont}
\usepackage{mathtools}
\usepackage{cite}
\usepackage{tikz-cd}
\usepackage{amssymb}
\usepackage{comment}
\usepackage{epigraph}
\hypersetup{
  colorlinks,
  citecolor=Blue,
  linkcolor=Black,
  urlcolor=Green}
  
\theoremstyle{plain}

\newtheorem{theorem}{Theorem}

\newtheorem{dummy}{anything}[section]

\newtheorem{lemma}[dummy]{Lemma}

\newtheorem{proposition}[dummy]{Proposition}

\theoremstyle{definition}
\newtheorem{definition}[dummy]{Definition}

\newtheorem{example}[dummy]{Example}

\newtheorem{remark}[dummy]{Remark}

\theoremstyle{remark}
\newcommand{\C}{\mathbb{C}}

\def\C{\mathbb{C}}

\title{Pseudo-holomorphic dynamics in the restricted three-body problem} 

\author{Agustin Moreno}

\address[A.\ Moreno]{School of Mathematics \\ Institute for Advanced Study \\ Princeton NJ  \\ USA}

\email{agustin.moreno2191@gmail.com}

\date{}

\begin{document}

\maketitle

\begin{abstract} In this article, we identify the 5-dimensional analogue of the finite energy foliations introduced by Hofer--Wysocki--Zehnder for the study of 3-dimensional Reeb flows, and show that these exist for the spatial circular restricted three-body problem (SCR3BP) whenever the planar dynamics is convex. We introduce the notion of a fiberwise-recurrent point, which may be thought of as a symplectic version of the leafwise intersections introduced by Moser, and show that they exist in abundance for a perturbative regime in the SCR3BP. We then use this foliation to induce a Reeb flow on the standard 3-sphere, via the use of pseudo-holomorphic curves, to be understood as the best approximation of the given dynamics that preserves the foliation. We discuss examples, further geometric structures, and speculate on possible applications.
\end{abstract}

\tableofcontents

\section{Introduction} The aim of this article is to delve into the well-known restricted three-body problem, concerning the motion of a small mass under the influence of two heavy masses, from the modern perspective of symplectic geometry.  We will focus on the \emph{spatial} problem, where the small mass moves in three-space, as opposed to the \emph{planar} one, where it moves in the plane. The main tools we will employ are the techniques from holomorphic curve theory and contact geometry. We will focus our attention on the energy level sets for energy below the first critical value, near the heavy masses. These carry contact structures by \cite{AFvKP}. The results in this paper fit into the scope of the \emph{symplectic dynamics} as proposed in \cite{BH}. We will:

\begin{enumerate}
    \item Identify the ``correct'' generalization in dimension $5$ of the finite energy foliations introduced by Hofer--Wysocki--Zehnder in dimension $3$ \cite{HWZ98}, in the sense that they indeed appear in the energy level set of the spatial (circular) restricted three-body problem, near each of the heavy masses, whenever the planar problem is convex. Each leaf is compatible with the dynamics, i.e.\ the dynamics induces a symplectic form which is positive on each annulus of the foliation. All such annuli have the same boundary (see Theorem \ref{thm:foliation}). In the case of the 3BP, this is the Hopf link corresponding to the direct/retrograde planar orbits;
    \item Identify the leaf space of the foliation (the moduli space of curves), i.e.\ it is the three-sphere $S^3$;
    
    \item Induce a geometric structure on the $3$-dimensional leaf space, in the form of a contact structure, which is obtained by an averaging procedure of the contact structure on the $5$-dimensional level set (see Theorem \ref{thm:contactstructure});
    
    \item Induce a dynamics on this leaf space, also induced by an averaging procedure of the original one, as well as show that every dynamics of the particular type obtained can be lifted (see Theorem \ref{thm:lifting}). This can be thought of as a statement that the \emph{complexity} of the dynamics in the space where the spatial problem occurs is higher than that of the space in which the planar problem occurs;
    
    \item Introduce a notion of \emph{fiberwise-recurrent} point, as a point whose orbit comes back to the leaf containing it. This can be thought of as a symplectic version of the notion of leafwise intersection of Moser, as the annuli are symplectic. 
    
    \item Prove the existence of infinitely many such recurrent points, for the perturbative case where the system is near-integrable (see Theorem \ref{thm:application}).
    
    \item Discuss examples (e.g.\ the well-known Katok example \cite{K73}), and introduce the notion of a cone structure being weakly/strongly adapted to an open book (see Definition \ref{def:stronglyadapted} and Definition \ref{def:weakly}). The strong notion indeed appears in the leaf space of the foliation of the 3BP, and the weak one, in a curious example due to G\"odel of a spacetime which allows for ``time travel'', as discussed in Appendix \ref{app:cones}.
   
\end{enumerate}

\textbf{Extended introduction.} In \cite{MvK}, as a higher-dimensional generalization of the global surfaces of section considered by Poincar\'e in the planar problem, the author and Otto van Koert showed that, in the low-energy range and independently of the mass ratio, the SCR3BP admits global hypersurfaces of section adapted to the Reeb dynamics. The underlying manifold is $S^2\times S^3$, and the hypersurfaces are copies of $\mathbb{D}^*S^2$, the pages of an open book whose binding is the planar problem $\mathbb{R}P^3$. From Appendix A in \cite{MvK}, we gather that in the integrable limit case of the (rotating) Kepler problem, the annuli fibers of the standard Lefschetz fibration on the page $P=\mathbb{D}^*S^2$ are invariant under the return map, which acts as a classical integrable twist map. Moreover, the moduli space of such fibers is naturally a copy of $S^3$, endowed with its trivial open book. 

\begin{figure}
    \centering
    \includegraphics[width=0.65\linewidth]{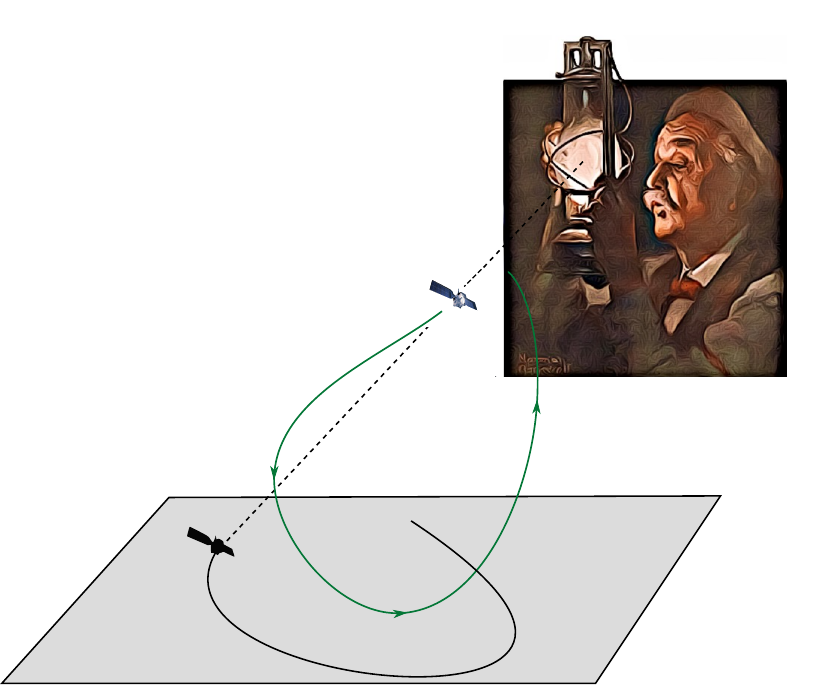}
    \caption{\textbf{Philosophy:} To shed some light on a complicated higher-dimensional problem, try first to look at the shadow that your lantern is producing!}
    \label{fig:lantern}
\end{figure}

In this article, we will generalize this geometric situation to a non-perturbative setting and for $5$-dimensional contact manifolds admitting an \emph{iterated planar} (IP) structure (as is the case of $S^2\times S^3$), with a view towards the SCR3BP. We first construct a version of a finite energy foliation of the $5$-dimensional manifold, which should be understood as the correct higher-dimensional analogue of the holomorphic open books as considered by Hofer-Wysocki-Zehnder \cite{HWZ98}, in order to study $3$-dimensional Reeb flows. This foliation may be thought of as an $S^1$-family of Lefschetz fibrations, one on each page of the open book, all inducing the same open book in the $3$-dimensional binding. We show that such a foliation always exists for the SCR3BP, whenever the planar dynamics is dynamically convex.

We then take a further step. To a Reeb dynamics on a contact $5$-fold adapted to an IP open book (i.e.\ its page admits a Lefschetz fibration with planar fibers), whose binding also carries an adapted open book, we associate a Reeb dynamics on a moduli space of holomorphic curves (a copy of $(S^3,\xi_{std})$), which is in some sense a ``holomorphic avatar'' or ``shadow'' of the original dynamics. Via this construction, the Kepler problem and its rotating version both correspond to the Hopf flow. When combined with \cite[Thm.\ 1.18]{HSW} (cf.\ \cite{HWZ98}; see also \cite{AFFHvK}), which guarantees the existence of adapted open books with annuli-like pages for the planar problem whenever the planar dynamics is dynamically convex, we obtain a holomorphic shadow for the SCR3BP for mass-ratio/energy pair in the convexity range (i.e.\ when the Levi-Civita regularization is strictly convex). We remark that convexity is not strictly needed, since all one needs for this construction is an adapted open book for the planar problem, and so dynamical convexity up to large action would also suffice (e.g.\ by perturbing the rotating Kepler problem, cf.\ \cite{AFFvK}).

The general direction is then to extract dynamical information for the $5$-fold, from information on its shadow; see Figure \ref{fig:lantern}. With this motivation in mind, we will then focus on global properties of this correspondence. In particular, we show that every Reeb dynamics on $S^3$ adapted to a concrete trivial open book arises as the holomorphic shadow of some Reeb dynamics on any given IP contact $5$-fold. We further study an non-trivial example due to Katok, where we explicitly relate the holomorphic shadow to suitable irrational ellipsoids. We also obtain dynamical information for the SCR3BP for the perturbative case where the mass ratio is sufficiently small, via the notion of a fiber-wise recurrent point, introduced below, which is a symplectic version of the well-known leafwise intersections introduced by Moser \cite{Moser} for the case of the isotropic foliation of a coisotropic submanifold.

\smallskip

\textbf{Setup.} Consider a \emph{concrete} open book decomposition $\pi: M \backslash B \rightarrow S^1$ on a contact $5$-manifold $(M,\xi)$. By definition, this is a fibration which coincides with the angle coordinate on a choice of collar neighborhood $B\times \mathbb{D}^2$ for the codimension-2 closed submanifold $B\subset M$ (the \emph{binding}). We assume that it supports $\xi$ in the sense of Giroux. This means that there is a contact form $\alpha$ for $\xi$, a \emph{Giroux form}, such that $\alpha\vert_B$ is contact, and $d\alpha$ is positively symplectic on the fibers of $\pi$; equivalently, the Reeb flow of $\alpha$ has $B$ as an invariant subset, and it is positively transverse to each fiber (see e.g.\cite{MvK} for more precise definitions). We denote the \emph{$\theta$-page} by $P_\theta=\overline{\pi^{-1}}(\theta)$ for $\theta \in S^1$, and we also use the abstract notation $M=\mathbf{OB}(P,\phi)$, where $P$ is the abstract page (the closure of the typical fiber of $\pi$) with $\partial P=B$, and $\phi$ is the symplectic monodromy. We assume that $P$ (abstractly) admits the structure of a $4$-dimensional Lefschetz fibration over $\mathbb{D}^2$ whose fibers are surfaces of genus zero and perhaps several boundary components. We abstractly write $P=\mathbf{LF}(F,\phi_F)$, where $\phi_F$ is the monodromy of the Lefschetz fibration on $P$ (necessarily a product of positive Dehn twists on the genus zero surface $F$). 

Following \cite{Acu}, we will refer to the open book on $M$ as an \emph{iterated planar} (IP) open book decomposition, and the contact manifold $M$ as iterated planar. As observed in \cite[Lemma 4.1]{AEO}, a contact $5$–manifold is iterated planar if and only if it admits an open book decomposition supporting the contact structure, whose binding is planar (i.e.\ admits a $3$-dimensional supporting open book whose pages have genus zero). In fact, we have $B=\mathbf{OB}(F,\phi_F)$. 

We wish to adapt the underlying planar structure to a \emph{given} Reeb dynamics on $M$ (and hence the need to work with concrete open books, rather than the abstract version). We then assume that the concrete open book on $M$ is adapted to the Reeb dynamics of a \emph{fixed} contact form $\alpha$, i.e.\ $\alpha$ is a Giroux form for the open book (whose dynamics we wish to study). In particular, $\omega_\theta:=d\alpha\vert_{P_\theta}$ is a symplectic form on $P_\theta$ for each $\theta \in S^1$. Therefore $(P_\theta,\omega_\theta)$ is an ideal Liouville filling of the binding $(B, \xi_B=\ker \alpha_B)$, where $\alpha_B=\alpha\vert_B$, for each $\theta$. We will further assume that we have a \emph{concrete} planar open book on the $3$-manifold $B=\mathbf{OB}(F,\phi_F)$, which is adapted to the Reeb dynamics of $\alpha_B$ and where $\phi_F$ is a product of positive Dehn twists in the genus zero surface $F$. We will denote $L=\partial F$, which is a link in $B$ (the binding of the open book for $B$, and Reeb orbits for $\alpha_B$). Given the above situation, we will say that the Giroux form $\alpha$ is an \emph{IP Giroux form}.

This is precisely the situation in the SCR3BP whenever the planar dynamics is convex/dynamically convex, as follows from \cite[Thm.\ 1.18]{HSW}, combined with Theorem 1 in \cite{MvK}.

\medskip

\textbf{Statement of results.} Our first result provides a foliation of $M$ via symplectic surfaces (with respect to $d\alpha$), which also foliate each page of the given IP open book, and which are fibers of an $S^1$-family of Lefschetz fibrations. 

    \begin{figure}
        \centering
        \includegraphics[width=0.5 \linewidth]{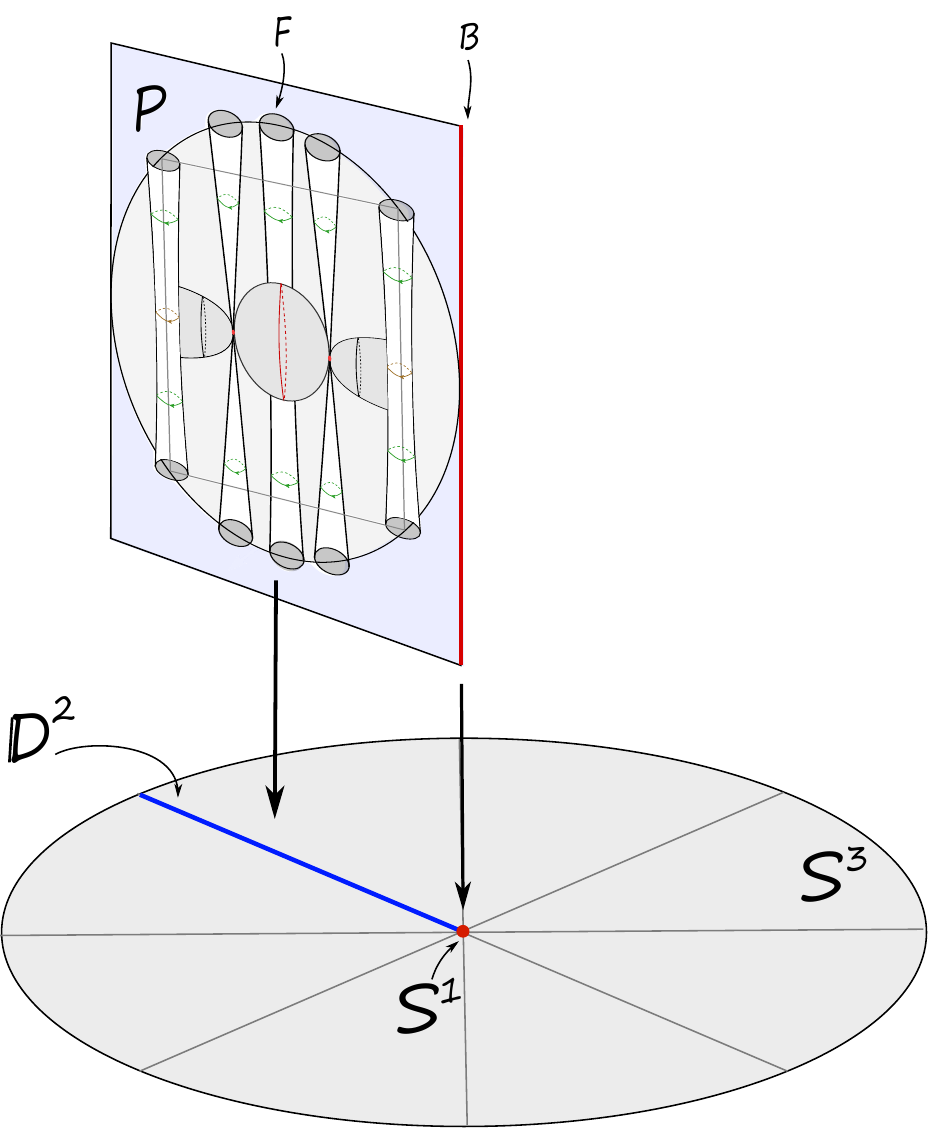}
        \caption{The moduli space of curves is a copy of $S^3=\mathbf{OB}(\mathbb{D}^2,\mathds{1})$.}
        \label{fig:MODULI}
    \end{figure}

\begin{theorem}[IP foliation]\label{thm:foliation}
There is a foliation $\overline{\mathcal{M}}^q$ of $M\backslash L$, consisting of immersed symplectic surfaces for $d\alpha$ whose boundary is $L$. Away from $B$, its elements are arranged as fibers of Lefschetz fibrations $\pi_\theta: P_\theta \rightarrow \mathbb{D}^2_\theta$, $\theta \in S^1$, all of which induce the same fixed concrete open book at $B$. The $\pi_\theta$ are all generic, i.e.\ each fiber contains at most a single critical point. We have $\overline{\mathcal{M}}^q\cong S^3=\mathbf{OB}(\mathbb D^2,\mathds 1)$, endowed with a trivial open book $\theta_\mathcal{M}: \overline{\mathcal{M}}^q\backslash \mathcal{M}_B^q\rightarrow S^1$ whose $\theta$-page is identified with $\mathbb{D}^2_\theta$, and its binding is $\mathcal{M}_B^q\cong S^1$, the family of pages of the open book at $B$. \end{theorem}

See Figure \ref{fig:MODULI} for a sketch of the geometric picture. We call a foliation as above, an \emph{IP foliation}.

\medskip

\textbf{Contact and symplectic structures on moduli.} The leaf-space of the above foliation comes endowed with extra structure, naturally induced from the structure on the $5$-fold. From Theorem \ref{thm:foliation}, we may view $\overline{\mathcal{M}}^q=\mathbf{OB}(\mathbb{D}^2,\mathds{1})\cong S^3$, and $\overline{\mathcal{M}}=\mathbb{R}\times \overline{\mathcal{M}}^q$. 

\begin{theorem}[contact and symplectic structures on moduli]\label{thm:contactstructure} The moduli space $\overline{\mathcal{M}}^q$ carries a natural contact structure $\xi_\mathcal{M}$ which is supported by the trivial open book on $S^3$ (and hence it is isotopic to the standard contact structure $\xi_{std})$. Moreover, the symplectization form on $\mathbb{R}\times M$ associated to any Giroux form $\alpha$ on $M$ induces a tautological symplectic form on $\overline{\mathcal{M}}$ by leaf-wise integration, which is naturally the symplectization of a contact form $\alpha_\mathcal{M}$ for $\xi_{\mathcal{M}}$, whose Reeb flow is adapted to the trivial open book on $\overline{\mathcal{M}}^q$.

\end{theorem}

Roughly speaking, the contact distribution $\xi_\mathcal{M}$ is induced by a $2$-dimensional subdistribution of $\xi=\ker \alpha$, which along $P_\theta$ is naturally isomorphic to a symplectic connection for the Lefschetz fibration $\pi_\theta$.

\begin{remark}[Holomorphic curves in moduli]
The above construction opens up the possibility of studying punctured holomorphic curves in the moduli space $\overline{\mathcal{M}}=\mathbb{R}\times \overline{\mathcal{M}}^q$. These correspond to $4$-dimensional $J$-invariant and asymptotically cylindrical hypersurfaces in $\mathbb{R}\times M$ (in the sense of \cite{MS}). For example, the holomorphic open book construction (see below) on the trivial open book for $\overline{\mathcal{M}}^q$ recovers the codimension-$2$ foliation $\mathcal{F}$ on $\mathbb{R}\times M$, used below for the proof of the above theorem.
\end{remark}

\textbf{The holomorphic shadow.} We define the \emph{(pseudo-)holomorphic shadow} of the Reeb dynamics of $\alpha$ on $M$ to be the Reeb dynamics of the associated contact form $\alpha_\mathcal{M}$ on $S^3$, provided by Theorem \ref{thm:contactstructure}. The flow of $\alpha_{\mathcal{M}}$ can be viewed as a flow $\phi_t^{M;\mathcal{M}}$ on $M\backslash L$ which leaves the holomorphic foliation $\overline{\mathcal{M}}^q$ invariant (i.e.\ it maps holomorphic curves to holomorphic curves). It is the ``best approximation'' of the Reeb flow of $\alpha$ with this property, as its generating vector field is obtained by reparametrizing and projecting the original Reeb vector field to the tangent space of $\overline{\mathcal{M}}^q$, via a suitable $L^2$-orthogonal projection. It may also be viewed as a Reeb flow $\phi_t^{S^3;\mathcal{M}}$ on $S^3$, related to the one on $M$ via a semi-conjugation
\begin{center}\label{diag:semiconj}
    \begin{tikzcd}
        M\backslash L\arrow[r, "\phi_t^{M;\mathcal{M}}"] \arrow[d, "\pi^q"]&
        M\backslash L \arrow[d, "\pi^q"]\\
        S^3 \arrow[r, "\phi_t^{S^3;\mathcal{M}}"]& S^3 \\
    \end{tikzcd}
\end{center}
where $\pi^q$ is the projection to the leaf-space $\overline{\mathcal{M}}^q\cong S^3$. 

\smallskip

\begin{figure}
    \centering
    \includegraphics[width=0.7 \linewidth]{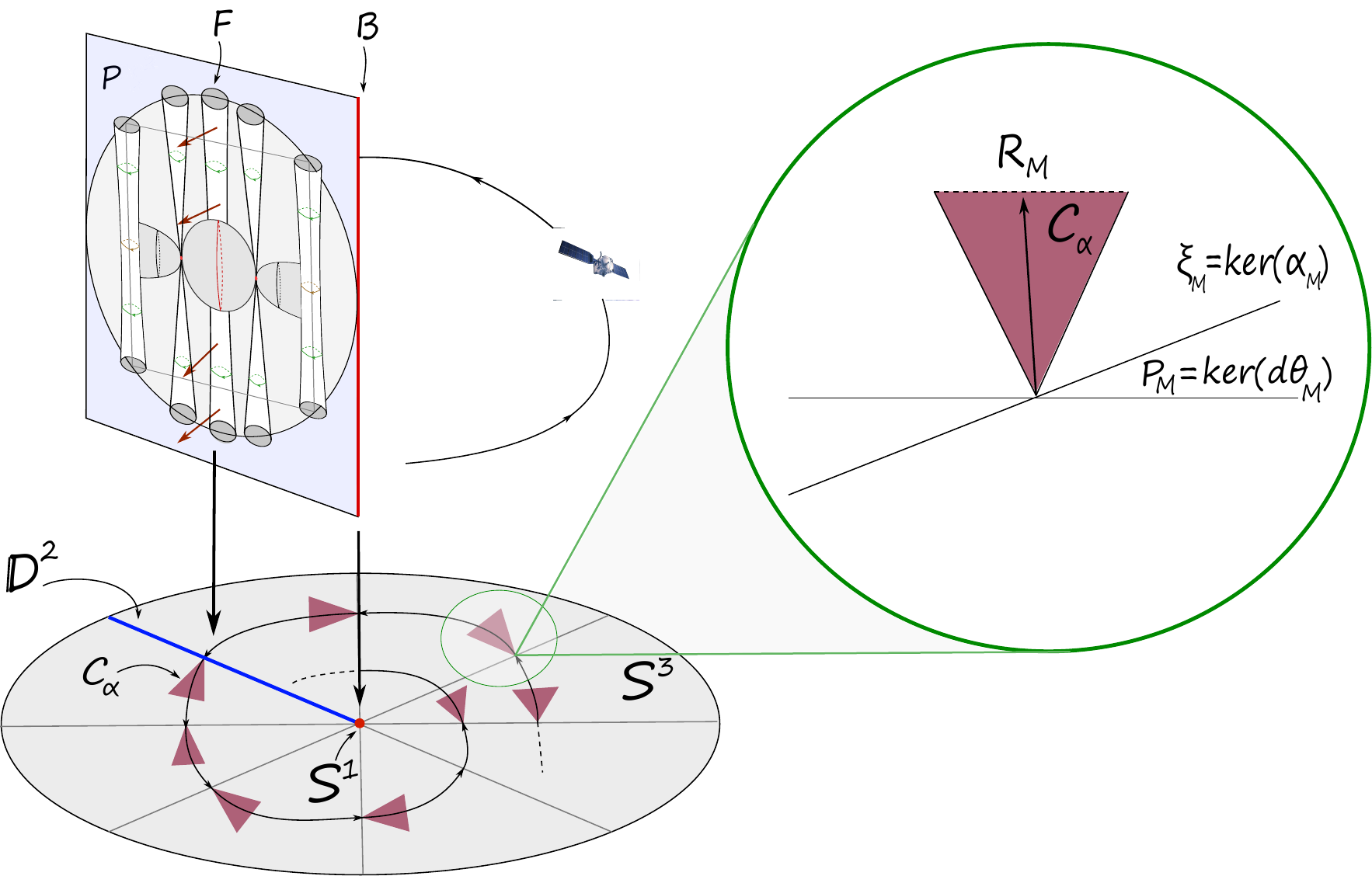}
    \caption{The shadowing cone is $C_\alpha=d\pi(\ker d\alpha)$. Orbits of $\alpha$ project to orbits of the cone, which are transverse to $\xi_\mathcal{M}$ and to every page. The Reeb vector field $R_\mathcal{M}$ spans the average direction of $C_\alpha$.}
    \label{fig:cone}
\end{figure}  

\textbf{Shadowing cone.} While the contact form $\alpha_\mathcal{M}$ is obtained from $\alpha$ by an averaging procedure that modifies the original dynamics by forcing it to preserve a foliation, another piece of the underlying geometric structure which encodes more reliable dynamical information is the \emph{shadowing cone} $C_\alpha=d\pi^q(\ker d\alpha)$. By construction, orbits of $\alpha$ project under $\pi^q$ to orbits of $C_\alpha$. Moreover, this cone is adapted to the open book on $S^3$, in the sense that along $\mathcal{M}_B^q$ we have $C_\alpha\vert_{\mathcal{M}_B^q}=T\mathcal{M}^q_B$, and such that $d\theta_\mathcal{M}$ and $\alpha_\mathcal{M}$ are positive on $C_\alpha$ away from $\mathcal{M}_B^q$; see Definition \ref{def:stronglyadapted} below. As $\alpha_\mathcal{M}$ is obtained by fiber-wise integration with respect to $\pi^q$, the Reeb vector field $R_\mathcal{M}$ of $\alpha_\mathcal{M}$ spans the average direction in $C_\alpha$. See Figure \ref{fig:cone}. One can therefore think of the holomorphic shadow as the ``guiding direction'' of the cone. For future use, we encode the properties of this cone in the following general definition:

\begin{figure}
    \centering
    \includegraphics[width=0.5 \linewidth]{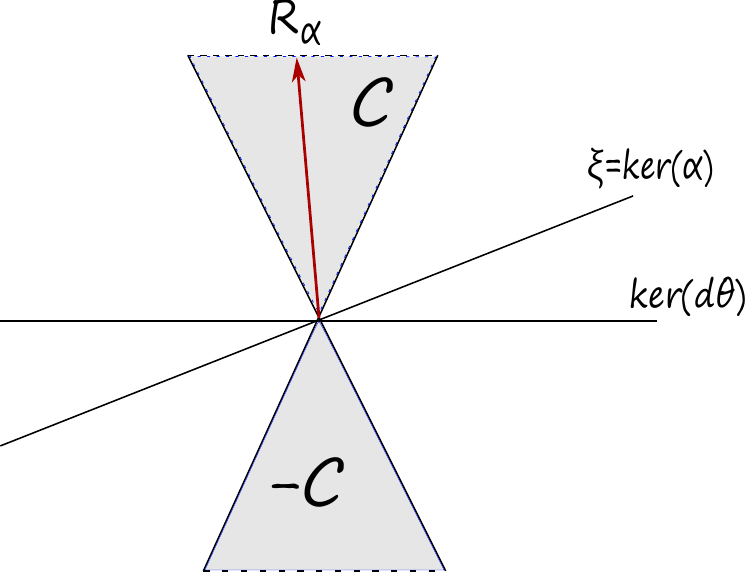}
    \caption{A cone structure adapted to an open book.}
    \label{fig:adapted}
\end{figure}

\begin{definition}\label{def:stronglyadapted}
Consider an everywhere non-trivial cone structure $C$ on a manifold $M$, where $M$ is endowed with an open book $\theta: M\backslash B\rightarrow S^1$. We say that $C$ is strongly adapted (or simply adapted) to $\theta$ if
\begin{enumerate}
    \item[(1)] $C\vert_B\subset TB$;
    \item[(2)] $d\theta$ is a section of $C\vert_{M\backslash B}$; 
\end{enumerate}
and if there exists a Giroux form $\alpha$ for the open book such that
\begin{enumerate}
    \item [(3)] $\alpha$ is a section for $C$;
    \item [(4)] The Reeb vector field $R_\alpha$ belongs to $C$.
\end{enumerate}
Here, a section for $C$ is a $1$-form which is strictly positive on non-zero vectors of $C$. See Figure \ref{fig:adapted}.
\end{definition}

In the case where $M$ is oriented, each page is co-oriented by the open book and hence inherits an orientation; then $B$ does also, as the boundary of each page. Then, if $M$ is $3$-dimensional and oriented, the first condition means that $C\vert_B=TB^+$, where $TB^+$ is defined as those vectors in $TB$ which point in the non-negative orientation of $B$. Therefore, each positively parametrized circle in the binding is an orbit for $C$.

\textbf{Global properties.} We will now focus on the global properties of the correspondence $\alpha \mapsto \alpha_\mathcal{M}$. For $F$ a genus zero surface, let $\mathbf{Reeb}(F,\phi_F)$ denote the collection of contact forms whose flow is adapted to some concrete planar open book $\pi_B: B\backslash L \rightarrow S^1$ on a given $3$-manifold $B$, of abstract form $B=\mathbf{OB}(F,\phi_F)$. Iteratively, we define $\mathbf{Reeb}(\mathbf{LF}(F,\phi_F),\phi)$ to be the collection of contact forms with flow adapted to some concrete IP open book $\pi_M: M\backslash B \rightarrow S^1$ on a $5$-manifold $M$, of abstract form $M=\mathbf{OB}(\mathbf{LF}(F,\phi_F),\phi)$, whose restriction to the binding $B=\mathbf{OB}(F,\phi_F)$ belongs to $\mathbf{Reeb}(F,\phi_F)$. We call elements in $\mathbf{Reeb}(\mathbf{LF}(F,\phi_F),\phi)$ \emph{IP contact forms}, or \emph{IP Giroux forms}. We may topologize both spaces $\mathbf{Reeb}(F,\phi_F)$ and $\mathbf{Reeb}(\mathbf{LF}(F,\phi_F),\phi)$ with the $C^\infty$-topology, for which both are infinite-dimensional CW complexes. 

We then have a map
$$
\mathbf{HS}:\mathbf{Reeb}(\mathbf{LF}(F,\phi_F),\phi)\rightarrow \mathbf{Reeb}(\mathbb{D}^2,\mathds{1}),
$$
$$
\alpha \mapsto \alpha_\mathcal{M},
$$
given by taking the holomorphic shadow. 

\begin{remark} The contact form $\alpha_\mathcal{M}$ map depends on a choice of almost complex structure $J$ on $\mathbb{R}\times M$, compatible with $d\alpha$ along $\ker \alpha$, and making $\mathbb{R}\times B$ a $J$-invariant submanifold (the $\mathbb R$-direction, however, has to be mapped by $J$ to the Reeb direction of a suitable SHS deformation of $\ker \alpha$, as explained below). The space of such choices, which we denote by $\mathcal{J}_{\alpha,B}$, is contractible for every $\alpha$, and hence $\mathbf{HS}$ depends on this choice only up to homotopy. Therefore the domain of $\mathbf{HS}$, stictly speaking, is
$$
\mathcal{D}_{\mathbf{HS}}:=\{(\alpha,J): \alpha \in \mathbf{Reeb}(\mathbf{LF}(F,\phi_F),\phi),\; J \in \mathcal{J}_{\alpha,B}\}. 
$$
We will write $\mathbf{HS}(\alpha,J)$ for the holomorphic shadow of $\alpha$ whenever we wish to emphasize this dependence.
\end{remark}

We will refer to $\mathbf{HS}^{-1}(\alpha_{std})$ as the \emph{integrable fiber} (note that this consists of pairs $(\alpha,J),$ as explained in the above remark).

\begin{theorem}[Reeb flow lifting theorem]\label{thm:lifting} $\mathbf{HS}$ is surjective.
\end{theorem}

In other words, we may lift any Reeb flow on $S^3$ adapted to the trivial open book, as the holomorphic shadow of the Reeb flow of an IP Giroux form adapted to \emph{any} choice of concrete IP contact $5$-fold (for some auxiliary choice of $J$). The map $\mathbf{HS}$ is clearly not in general injective (i.e.\ the shadow does not see vertical dynamical information), as e.g.\ Example \ref{ex:rotKepler} shows. The above theorem says that Reeb dynamics on \emph{any} IP contact $5$-fold is at least as complex as Reeb dynamics on the standard contact $3$-sphere (and so provides a concrete measure of the high complexity of, say, the spatial restricted three-body problem, i.e.\ it is at least as complicated as the planar problem). We point out that higher-dimensional Reeb flows encode the complexity of all flows on arbitrary compact manifolds (i.e.\ they are \emph{universal}) \cite{CMPP}. See also the related discussion at the end of the article on topological entropy.

\medskip

\textbf{Dynamical Applications.} We wish to apply the above results to the SCR3BP. We first introduce the following general notion. Consider an IP $5$-fold $M$ with an IP Reeb dynamics, endowed with an IP holomorphic foliation $\mathcal{M}$ as in Theorem \ref{thm:foliation}. Fix a page $P$ in the IP open book of $M$, and consider the associated Poincar\'e return map $f: \mbox{int}(P) \rightarrow \mbox{int}(P)$. A (spatial) point $x\in \mbox{int}(P)$ is said to be \emph{leaf-wise} (or \emph{fiber-wise}) \emph{$k$-recurrent} with respect to $\mathcal{M}$ if $f^k(x) \in \mathcal{M}_x$, where $\mathcal{M}_x$ is the leaf of $\mathcal{M}$ containing $x$, and $k\geq 1$. This means that $f^k(\mbox{int}(\mathcal{M}_x))\cap \mbox{int}(\mathcal{M}_x)\neq \emptyset$. This is, roughly speaking, a symplectic version of the notion of \emph{leaf-wise intersection} introduced by Moser \cite{Moser} for the case of the isotropic foliation of a coisotropic submanifold. 

\begin{figure}
    \centering
    \includegraphics[width=1 \linewidth] {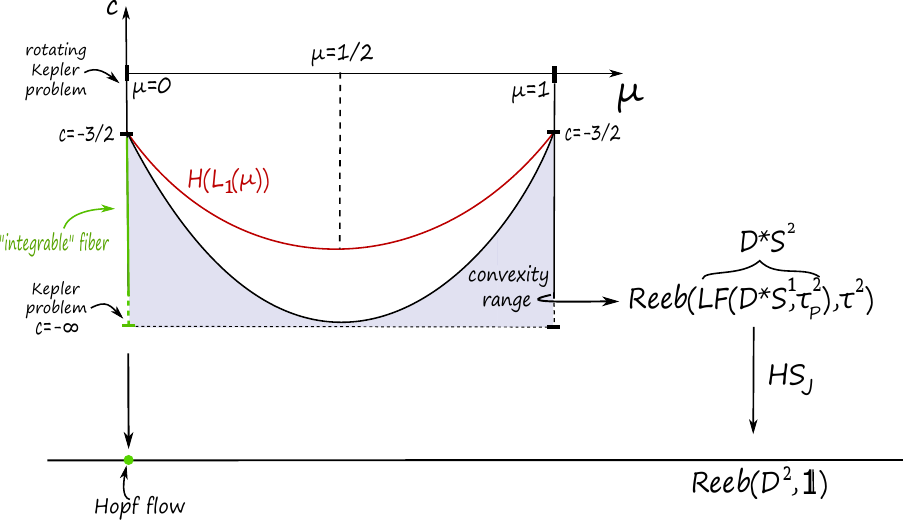}
    \caption{An abstract sketch of the convexity range in the SCR3BP (shaded), for which the holomorphic shadow is well-defined. Here, $c$ is the Jacobi constant, $\mu$ is the mass ratio, $\tau_P$ is the Dehn twist on $\mathbb{D}^*S^1$, and $\tau$ is the Dehn-Seidel twist in $\mathbb{D}^*S^2$. We should disclaim that the above is not a plot; the convexity range is not yet fully understood, although it contains (perhaps strictly) a region which qualitatively looks like the above, cf.\ \cite{AFFHvK,AFFvK}.}
    \label{fig:my_label}
\end{figure}

In the integrable case of the rotating Kepler problem, where the mass ratio $\mu=0$, the holomorphic foliation provided by Theorem \ref{thm:foliation} can be obtained by an explicit construction on $S^*S^3$ of a Lefschetz fibration on $T^*S^3$ with two singular fibers, such that the nodal singularities are fixed points of the symplectic return map; see \cite[Appendix A]{MvK}. Denote this ``integrable'' holomorphic foliation on $S^*S^3$ by $\mathcal{M}_{int}$. Since the return map for $\mu=0$ preserves fibers, every point is leaf-wise $1$-recurrent with respect to $\mathcal{M}_{int}$ (and the shadowing cone is the positive span of the Hopf flow). If the mass ratio is sufficiently small, then the leaves of $\mathcal{F}_{int}$ will still be symplectic with respect to $d\alpha$, where $\alpha$ is the corresponding perturbed contact form on the unit cotangent bundle $S^*S^3$. 

In Section \ref{sec:tomographies}, we introduce the notion of the \emph{transverse shadow}, which records which holomorphic curve is intersected by the orbit of each point, and therefore consists of paths in $S^3$ which are positively transverse to the standard contact structure $\xi_\mathcal{M}$ and to the pages of the open book. The transverse shadow is just the collection of orbits of the shadowing cone which come from orbits of the original dynamics. By coherently packaging these transverse paths using the notion of a \emph{symplectic tomography} also introduced in said section, and appealing to the classical Brouwer's translation theorem on the open disk, we will obtain the following perturbative result: 

\begin{theorem}\label{thm:application}
In the SCR3BP, for any choice of page $P$ in the open book of \cite{MvK}, for any fixed choice of $k\geq 1$, for sufficiently small $\mu \gtrsim 0$ (depending on $k$), for energy $c$ below the first critical value $H(L_1(\mu))$, along the bounded component of the Hill region near the Earth, and for every $l\leq k$, there exist infinitely many points in $\mbox{int}(P)$ which are leaf-wise $l$-recurrent with respect to $\mathcal{M}_{int}$.
\end{theorem}

In simpler words, there are plenty of leaf-wise recurrent points near the Earth, when the Moon is small.

\begin{remark}
The same conclusion holds for arbitrary $\mu \in [0,1]$, but sufficiently negative $c\ll 0$ (depending on $\mu$ and $k$).
\end{remark}

In fact, the conclusion of the Theorem \ref{thm:application} holds whenever the relevant return map is sufficiently close to a return map which preserves the leaves of the holomorphic foliation of Theorem \ref{thm:foliation} (i.e.\ which coincides with its holomorphic shadow on $M$). It may then be interpreted as a symplectic version of the main theorem in \cite{Moser}, for two-dimensional symplectic leaves. The advantage of considering the integrable foliation (in terms of applications) is that it can be qualitatively understood, as in \cite[Appendix A]{MvK}; the above conclusion also holds for the corresponding perturbed foliation.

\medskip

\textbf{Acknowledgements.} The author is deeply grateful to Otto van Koert, for discussions and inputs on the material here presented, which stems from their joint work. To Urs Frauenfelder, who encouraged the start of this long project and followed its developments closely. To Helmut Hofer, Kai Cieliebak, Umberto Hrnyewicz, Alejandro Passeggi, for further discussions and inputs throughout different stages of this program. To the anonymous referee, for the careful summary of the results here presented. The author also acknowledges a Research Fellowship funded by the Mittag-Leffler Institute in Djursholm, Sweden, where an earlier version of this manuscript was finalized. The author further acknowledges the support by the National Science Foundation under Grant No. DMS-1926686, as part of the funding for his membership at the Institute for Advanced Study in Princeton.

\section{Dynamics on moduli spaces}

Let $M=\mathbf{Reeb}(P,\phi)$ be an IP $5$-fold, $\alpha \in \mathbf{Reeb}(P,\phi)$, and $\xi=\ker \alpha$ with Reeb vector field $R_\alpha$. Let $\theta_B: B\backslash L\rightarrow S^1$ be a concrete open book on $B$ adapted to $\alpha_B=\alpha\vert_B$, and $\theta_M: M\backslash B\rightarrow S^1$ a concrete open book adapted to $\alpha$. We denote by $P_\theta=\theta_M^{-1}(\theta)$, $F_\theta=\theta_B^{-1}(\theta)$, the $\theta$-pages.

\smallskip

\textbf{Holomorphic foliations.} For the sake of brevity, we will streamline the arguments which have already appeared throughout the literature and give the appropiate references. We proceed to the construction for the proof of Theorem \ref{thm:foliation}. We may construct the fibers of a Lefschetz fibration $\pi_\theta: P_\theta \rightarrow \mathbb{D}^2_\theta$ on each page $P_\theta$, with regular fiber $F$, in such a way that these fibrations form a well-defined $S^1$-family in $M$; and so that such fibers foliate the contact manifold $M$, and are symplectic submanifolds in each page. We can do so as follows.

\smallskip

\textbf{Finite energy foliation over $B$.} First, we may construct an almost complex structure $J_B$ on the $4$-dimensional symplectization $\mathbb{R}\times B$, together with an $\mathbb{R}$-invariant and Fredholm regular $J_B$-holomorphic finite energy foliation $\mathcal{M}_B$ of $\mathbb{R}\times B$, by punctured holomorphic curves which are asymptotically cylindrical Liouville completions $\widehat F_\theta$ of $F_\theta$, asymptotic to $L$, so that $\widehat F_\theta$ projects to $F_\theta$ under the projection $\pi_B: \mathbb{R}\times B \rightarrow B$. Moreover, we have $\mathcal{M}_B^q:=\mathcal{M}_B/\mathbb R\cong S^1$. This construction of a ``holomorphic'' open book has appeared in several forms in the literature (e.g.\ \cite{W10,A11}), and therefore we will omit details. Nevertheless, there is a slight technicality: unlike e.g.\ in \cite{W10}, the symplectic form in $\mathbb{R}\times M$ is a priori \emph{given} as $\omega_B=d(e^t \alpha_B)$. In order to construct the holomorphic foliation, we need to deform $(\alpha_B,d\alpha_B)$ to a stable Hamiltonian structure (SHS) $\mathcal{H}_B=(\lambda_B, d\alpha_B)$ which is tangent to pages away from the binding, and the resulting $J_B$ will be compatible with the symplectization of $\mathcal{H}_B$, i.e.\ will map the $\mathbb{R}$-direction to the Reeb vector field of $\mathcal{H}_B$. The resulting holomorphic curves, while not strictly holomorphic for an almost complex structure compatible with $\alpha_B$, will \emph{still} be symplectic submanifolds for $\omega_B$, and this is what we ultimately care about for Theorem \ref{thm:foliation}. 

\smallskip

\begin{figure}
    \centering
    \includegraphics[width=0.5 \linewidth]{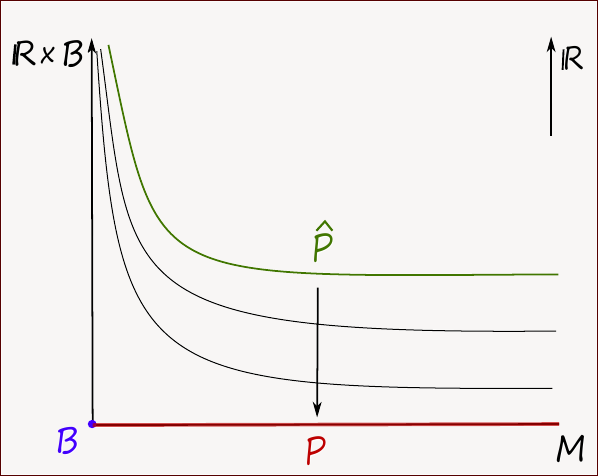}
    \caption{The codimension-2 holomorphic foliation $\mathcal{F}$.}
    \label{fig:foliation}
\end{figure}

\textbf{Codimension-$2$ foliation over $M$.} We may then do the same construction, but two-dimensions higher, as follows. Since $\xi$ is supported by the open book $\theta_M$, away from $B$ we have an isomorphism $TP_\theta\cong \xi\vert_{P_\theta}$, and on $B$ we have $\xi\vert_B=\xi_B \oplus \xi_B^\perp$ where $\xi_B=\ker \alpha_B$, and $\xi_B^\perp$ the $d\alpha$-symplectic complement inside $\xi$ of $\xi_B$. For each page, one then extends $J_B$ to an almost complex structure $J_\theta$ on $P_\theta$ (viewed as an ideal Liouville filling of $B$), compatible with $\omega_\theta=d\alpha\vert_{P_\theta}$, and generic as a $1$-parameter family. One can then deform $(\alpha,d\alpha)$ to a stable Hamiltonian structure $\mathcal{H}=(\lambda, d\alpha)$ with kernel $\xi_\mathcal{H}=\ker \lambda$ deforming the kernel of $\alpha$, by simply deforming the $1$-form, so that $\lambda$ interpolates between $\alpha$ near the binding and $d\theta_M$ away from it. We then have $\xi_\mathcal{H}\cong \xi$; $\xi_\mathcal{H}$ is tangent to $P_\theta$ away from a neighbourhood of $B$, so we still have the property $\xi\vert_{P_\theta}\cong \xi_\mathcal{H}\vert_{P_\theta} \cong TP_\theta$; and $\xi_\mathcal{H}=\xi$ near $B$. Moreover, the Reeb vector field of $\mathcal{H}$ still coincides with $R_\alpha$. One can then induce an $\mathcal{H}$-compatible almost complex structure $J$ on $\mathbb{R}\times M$ (i.e.\ it preserves $\xi_\mathcal{H}$ where it is $d\alpha$-compatible, and maps $\partial_t$ to $R_\mathcal{\alpha}$), with the extra properties that $J\vert_{\mathbb{R}\times B}=J_B$, the splitting $\xi\vert_B=\xi_B \oplus \xi_B^\perp$ is $J$-complex, and $J\vert_{\xi_\mathcal{H}\vert_{P_\theta}}$ corresponds to $J_\theta$ under the isomorphism $\xi_\mathcal{H}\vert_{P_\theta} \cong TP_\theta$. Note that we also have a natural $J^\xi$ which is compatible with $\alpha$, deforming $J$, for which the isomorphism $\xi\cong \xi_\mathcal{H}$ becomes complex, and $J^\xi(\partial_t)=R_\alpha$. One then uses the $\mathcal{H}$-compatible $J$ to make the open book $\theta_M$ holomorphic, i.e.\ we obtain a codimension-2 holomorphic foliation $\mathcal{L}$ of $\mathbb{R}\times M$ whose leaves are $\mathcal{L}_B=\mathbb{R}\times B$ and Liouville completions $\widehat P_\theta$ of the pages $P_\theta$. We have that $\widehat P_\theta$ projects to $P_\theta$ under the projection $\pi_M: \mathbb R\times M\rightarrow M$, and is asymptotically cylindrical to the $\mathbb{R}$-invariant hypersurface $\mathcal{L}_B$ in the sense of \cite{MS}. Moreover, the complex hypersurface $(\widehat P_\theta,J\vert_{\widehat P_\theta})$ is biholomorphic to $(P_\theta,J_\theta)$. See Figure \ref{fig:foliation} for a sketch. 

\smallskip

\textbf{Moduli space over $M$.} The moduli space $\mathcal{M}_B$ of holomorphic curves in $\mathcal{L}_B$, which is a copy of $\mathbb{R}\times S^1$, naturally extends to a moduli space $\mathcal{M}$ in $\mathbb{R}\times M$. Note that every curve in $\mathcal{M}$, while $J$-invariant but \emph{not} $J^\xi$-invariant, is $d\alpha$-symplectic. Moreover, an application of both $4$-dimensional \emph{and} higher-dimensional Siefring intersection theory as in \cite{MS,Mo} implies that: either the image of a curve in $\mathcal{M}$ lies completely in $\mathcal{L}_B$ and is a leaf of $\mathcal{M}_B$; or it lies completely in a leaf of $\mathcal{H}$ (and in particular is disjoint from $\mathcal{L}_B$ except ``at infinity'', intersecting it along $L$). For convenience of the reader, we outline the higher-dimensional version of the argument (as fully explained in \cite{MS}), which first restricts the behaviour of holomorphic curves to lie inside the leaves of the codimension-$2$ foliation, so that the $4$-dimensional techniques then apply. We remark that the higher-dimensional and $4$-dimensional argument is exactly the same.

\medskip

\textbf{Controlling curves: Siefring intersection theory.} As explained in \cite{MS}, there is a well-defined and homotopy-invariant \emph{Siefring intersection pairing} $u*H$ for a curve $u\in \mathcal{M}$ and a hypersurface $H\in \mathcal{F}$. We claim that this pairing vanishes: using homotopy invariance with respect to $H$ completely analogously as e.g.\ \cite[Lemma 4.9]{LW}, we have
$$
u*H=u*(\mathbb{R}\times B).
$$
From homotopy invariance with respect to $u$, it suffices to compute the above in the case where $u\in \mathcal{M}_B$ lies completely in $\mathbb{R}\times B$. Then we can appeal to the following intersection formula \cite{Sie}:
$$
u*(\mathbb{R}\times B)=\frac{1}{2}(\mu_N^\tau(u)-\#\Gamma(u)_{odd}),
$$
where $\mu_N^\tau(u)$ is the total normal Conley-Zehnder index of $u$ with respect to a trivialization $\tau$ of the symplectic normal bundle to $\mathbb{R}\times B$, and $\#\Gamma(u)_{odd}$ is the number of asymptotics of $u$ which have odd normal Conley-Zehnder index. In our case, each asymptotic of $u$ has normal index equal to $1$, and the claim that $u*H=0$ follows. Moreover, if $u$ did not lie in a leaf of $\mathcal{F}$, this pairing would be strictly positive, by the fact that interior intersections contribute positively, and intersections at infinity all vanish (the asymptotic winding numbers of $u$ all vanish, and these are the extremal ones). This would be a contradiction, and so we get the desired restrictions of elements in $\mathcal{M}$.

\begin{figure}
    \centering
    \includegraphics[width=0.5 \linewidth]{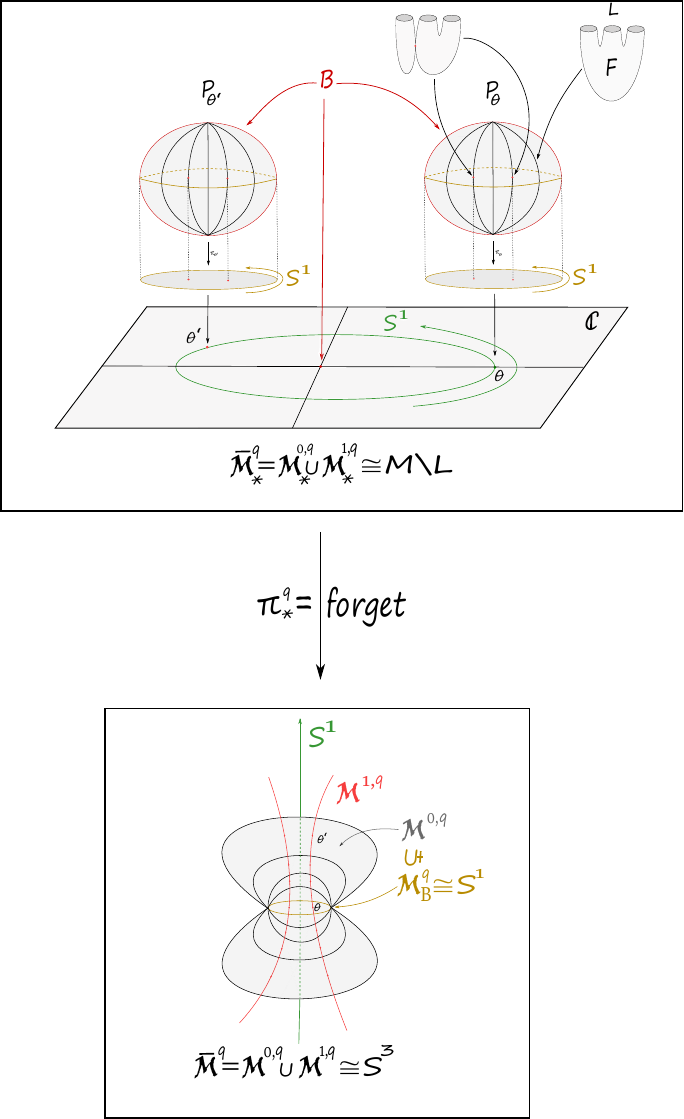}
    \caption{The forgetful map, and the stratification in the quotient moduli spaces. In the above picture, the boundary of each page is identified to each other.}
    \label{fig:forget}
\end{figure}

\smallskip

\textbf{SFT-Gromov compactification.} We may then consider the SFT-Gromov compactification $\overline{\mathcal{M}}$ of $\mathcal{M}$ by adding strata of nodal curves, as well as its pointed version $\overline{\mathcal{M}}_*$ by adding a single marked point in the domain of each curve, together with the resulting evaluation map $ev: \overline{\mathcal{M}}_* \rightarrow \mathbb{R}\times M$ which evaluates a curve at its marked point, and the forgetful map $\pi_*: \overline{\mathcal{M}}_*\rightarrow \overline{\mathcal{M}}$ which forgets the marked point. Moreover, one can apply Siefring intersection theory to the compactified moduli space too, and still conclude that the components of each element lie in the holomorphic hypersurfaces. Applying the results in \cite{W10b} (namely, Thm.\ 7 and its parametric version Thm.\ 8), one concludes that each $P_\theta$ admits the claimed generic Lefschetz fibration $\pi_\theta:P_\theta \rightarrow \mathbb{D}^2_\theta$, each inducing the original open book at $B$. Note that each such fibration is minimal (i.e.\ there are no contractible vanishing cycles), due to exactness of the symplectic form, and moreover every curve is immersed. In particular, non-nodal curves are embedded, and the nodal curves have embedded components, intersecting each other transversely; the same is true even after projecting via $\pi_M$.

This means that $\overline{\mathcal{M}}_*$ naturally admits a stratification $\overline{\mathcal{M}}_*=\mathcal{M}_*^0\bigsqcup \mathcal{M}_*^1$, where $\mathcal{M}_*^i$ consists of pointed curves with precisely $i$ nodes, and having closure $\overline{\mathcal{M}}_*^i=\bigcup_{j\leq i} \mathcal{M}_*^i$. In particular $\mathcal{M}_*^0$ is the top open stratum consisting of non-nodal curves, and $\mathcal{M}_*^1$ is closed. We similarly have a stratification for the unmarked moduli space $\overline{\mathcal{M}}=\mathcal{M}^0\bigsqcup \mathcal{M}^1$, where $\mathcal{M}^i=\pi_*(\mathcal{M}^i_*)$. The expected dimension of $\mathcal{M}_*^i$ is $6-2i$ (and that of $\mathcal{M}^i$ is $4-2i$), and from a similar analysis as carried out in \cite[sec.\ 4.7]{Mo} for curves lying in holomorphic hypersurfaces, one can show that each nodal strata $\mathcal{M}^i$ is Fredholm regular. 

\smallskip

\textbf{Structural diagrams.} The upshot is that we have a diagram of the form
\begin{center}
    \begin{tikzcd}
        \overline{\mathcal{M}}_*\arrow[r, "ev", "\cong"'] \arrow[d, "\pi_*"]&
        \mathbb{R}\times M\backslash \mathbb{R}\times L \arrow[dl, "\pi"]\\
        \overline{\mathcal{M}}\cong \mathbb{R}\times S^3& \\
    \end{tikzcd}
\end{center}
where $ev$ is a diffeomorphism since elements in $\overline{\mathcal{M}}$ together with trivial cylinders over orbits in $L$ foliate $\mathbb{R}\times M$, and so $\overline{\mathcal{M}}_*^q:=\overline{\mathcal{M}}_*/\mathbb{R}\cong M\backslash L$, and we denote $\pi=\pi_* \circ ev^{-1}$. We also used that $\overline{\mathcal{M}}\cong \mathbb{R}\times S^3$, and so $\overline{\mathcal{M}}^q:=\overline{\mathcal{M}}/\mathbb{R}\cong S^3$. Namely, $\overline{\mathcal{M}}^q$ is naturally equipped with the trivial open book $\overline{\mathcal{M}}^q=\mathbf{OB}(\mathbb{D}^2,\mathds{1})$, with $\theta$-page $P_\theta^\mathcal{M} \cong \mathbb{D}^2_\theta$ corresponding to the $2$-disk $\mathbb{D}^2_\theta$ at the base of the Lefschetz fibration $\pi_\theta$, via $\mathbb{D}_\theta^2=(\pi_\theta \circ ev^q \circ (\pi_*^q)^{-1})(P_\theta^\mathcal{M})$, where $ev^q$ and $\pi^q_*$ are the quotient maps induced by $ev$ and $\pi_*$, respectively. Its binding is identified with $\mathcal{M}_B^q:=\mathcal{M}_B/\mathbb{R}\cong S^1$. By quotienting out the $\mathbb{R}$-action, and denoting $\pi^q=\pi_*^q \circ (ev^q)^{-1}$, we obtain a similar diagram
\begin{center}
    \begin{tikzcd}
        \overline{\mathcal{M}}^q_*\arrow[r, "ev^q", "\cong"'] \arrow[d, "\pi^q_*"]&
        M\backslash L \arrow[dl, "\pi^q"] \\
        \overline{\mathcal{M}}^q\cong S^3& \\
    \end{tikzcd} 
\end{center}
See Figure \ref{fig:forget}, where we denote $\mathcal{M}_*^{q,i}=\mathcal{M}_*^{i}/\mathbb{R}$, $\mathcal{M}^{q,i}=\mathcal{M}^{i}/\mathbb{R}=\pi_*(\mathcal{M}_*^{q,i})$, for $i=0,1$.  We denote by $\theta_\mathcal{M}: \overline{\mathcal{M}}^q\backslash \mathcal{M}_B^q\rightarrow S^1$ the open book fibration on the quotient moduli space, which we call the \emph{shadowing open book}. This finishes the construction which proves Theorem \ref{thm:foliation}. 

\medskip

\textbf{Symplectic and contact forms on moduli spaces.} We now proceed to the construction underlying Theorem \ref{thm:contactstructure}. Let $$\mathcal{P}=\{\varphi \in C^\infty(\mathbb{R},(0,1)): \varphi^\prime>0\}$$ be the space of orientation-preserving diffeomorphisms between $\mathbb{R}$ and $(0,1)$. For each $\varphi \in \mathcal{P},$ we can induce a symplectic form $\omega_{\mathcal{M}}^\varphi$ on $\overline{\mathcal{M}}$ and a ($\varphi$-independent) contact form $\alpha_\mathcal{M}$ on $\overline{\mathcal{M}}^q$, such that $$(\overline{\mathcal{M}}\cong \mathbb{R}\times S^3,\omega^\varphi_{\mathcal{M}}=d(e^{\varphi(t)} \alpha_{\mathcal{M}}))$$ is the symplectization of the contact manifold $(\overline{\mathcal{M}}^q\cong S^3, \alpha_{\mathcal{M}}).$ Indeed, let $$\omega^\varphi:=d(e^{\varphi(t)}\alpha) \in \Omega^2(\mathbb{R}\times M)$$ for $\varphi \in \mathcal{P}$. Denote $\omega^\varphi_*:=ev^*\omega^\varphi \in \Omega^2(\overline{\mathcal{M}}_*)$, and $F_u=\pi_*^{-1}(u)$ the fiber over $u \in \overline{\mathcal{M}}$ (its domain). By construction, $\omega^\varphi_*$ is a symplectic form in $\overline{\mathcal{M}}_*$, which is symplectic on the fibers of $\pi_*$, i.e.\ $\omega_u^\varphi=\omega^\varphi_*\vert_{F_u}$ is an area form on the domain of $u$. Moreover, there is a tautological notion of a symplectic connection for $\pi_*$, so that every vector $v \in T_u \overline{\mathcal{M}}$ has a horizontal lift $\tilde{v} \in T_{(u,z)}\overline{\mathcal{M}}_*$, for every $(u,z) \in F_u$ (where $z$ lies in the domain of $u$). Indeed, a tangent vector field $v \in T_u \overline{\mathcal{M}}$ is simply a vector field along $F_u$, taking values in the (generalized; see below) normal bundle $N_u$ and lying in the kernel of the normal linearized CR-operator $\mathbf{D}_u^N$ at $u$. The horizontal lift $\widetilde{v}$ of $v$ at $(u,z)$ is simply $\widetilde{v}=v(z) \in N_u\vert_z$, i.e.\ the vector field itself. This also makes sense along nodal curves. Indeed, recall that the tangent space to a \emph{fixed} (smooth) stratum $\mathcal{M}^i$ of the moduli space along a nodal curve $u=(u_1,\dots,u_{i+1})$ consists of tuples $\eta=(\eta_1,\dots,\eta_{i+1})$ of normal sections $\eta_j$ along each component $u_j$, lying in the kernel of $\bigoplus_{j=1;EV}^{i+1}\mathbf{D}^N_{u_j}$, the fiber-product (under the linear evaluation map $EV$ at the nodes) of the normal linearized CR-operators at each $u_j$; so in particular each element in the tuple $\eta$ agrees at each nodal pair. Tangent vectors along a nodal curve which are not tangent to the fixed stratum containing it can also be thought of as vector fields along the curve (as follows from standard gluing analysis). So the above notion of a horizontal lift carries through immediately.

For $v,w \in T_u \overline{\mathcal{M}}$, let $\omega_*^\varphi(v,w):=\omega^\varphi_*(\widetilde{v},\widetilde{w})$ viewed as a function on $F_u$, and define $\omega^\varphi_{\mathcal{M}} \in \Omega^2(\overline{\mathcal{M}})$ via
$$
(\omega^\varphi_{\mathcal{M}})_{u}(v,w)=\int_{z\in F_u} \omega_*^\varphi(v(z),w(z))dz,
$$
where by simplicity we denote $dz=\omega_u^\varphi$. We remark that this same construction for the special case of closed and immersed curves (not including nodal degenerations) has been carried out in \cite{CKP}, although with a slightly different language. In general, generically and in dimension $6$, the locus of non-immersed curves consists of isolated curves with a single simple critical point (see Appendix A in \cite{W19} for the closed case). This, and the fact that there is still a well-defined normal bundle at critical points --the \emph{generalized} normal bundle \cite{W10c}--, implies that the above integral is well-defined. In our special setup, we have already seen that every curve is immersed, so this is not even an issue. While we moreover have punctures in our setup, the finiteness of the integral is ensured by our choice of a diffeomorphism $\varphi$. The fact that $\omega^\varphi_{\mathcal{M}}$ is symplectic follows by adapting the arguments in \cite{CKP} (except the nondegeneracy property: here $J$ is not necessarily integrable, and so the kernel of $\mathbf{D}^N_u$ is not necessarily a complex space. But we understand precisely what the kernel is, due to the foliation property and Fredholm regularity, and hence immediately see that non-degeneracy follows; see below for more details). The contact form on $\overline{\mathcal{M}}^q$ is obtained from the symplectic form in the obvious way, i.e.\ as 
$$\alpha_{\mathcal{M}}=(e^{\varphi}\varphi^\prime)^{-1}i_{\partial_t}\omega^\varphi_{\mathcal{M}} \in \Omega^1(\overline{\mathcal{M}}^q).
$$
By construction, this contact form is independent on the choice of $\varphi$. Indeed, we can give a more explicit description as follows. First note that $d\alpha\vert_{F_u}$ is a symplectic form on the interior of $F_u$ (which degenerates at the boundary since $\partial F_u=L$ is a collection of $\alpha$-orbits). The volume of $F_u$ with respect to $d\alpha\vert_{F_u}$ is, by Stokes' theorem, given by $$\mbox{vol}(F_u)=\int_{F_u}d\alpha\vert_{F_u}=\int_{L}\alpha=:A_\alpha>0,$$ the total $\alpha$-action of $L$, which is independent of $u$. We then have 
$$
(\alpha_\mathcal{M})_u(v)=\int_{z\in F_u}\alpha_z(v(z))dz,
$$
$$
(d\alpha_\mathcal{M})_u(v,w)=\int_{z\in F_u}d\alpha_z(v(z),w(z))dz,
$$
for $v,w\in T\overline{\mathcal{M}}^q$, where for simplicity we use the notation $dz$ to indicate that we integrate over $z$ with respect to the area form $d\alpha\vert_{F_u}$, and where we use the same letter $u$ to denote the projection of $u \in \mathcal{M}$ to $\mathcal{M}^q=\mathcal{M}\backslash \mathbb{R}$ (which is also embedded). The contact structure $\xi_{\mathcal{M}}=\ker \alpha_{\mathcal{M}}$ can be understood, over $u$, as the average of the contact structures $\xi_z$ for $z\in F_u$; see Lemma \ref{lemma:adapted} below. We call $\xi_\mathcal{M}$, the \emph{shadowing contact structure}. Moreover, it is supported by the trivial open book on $S^3$, and hence it is isotopic to the standard contact structure on $S^3$. This can be understood as follows. 

\smallskip

\begin{figure}
    \centering
    \includegraphics[width=0.5 \linewidth]{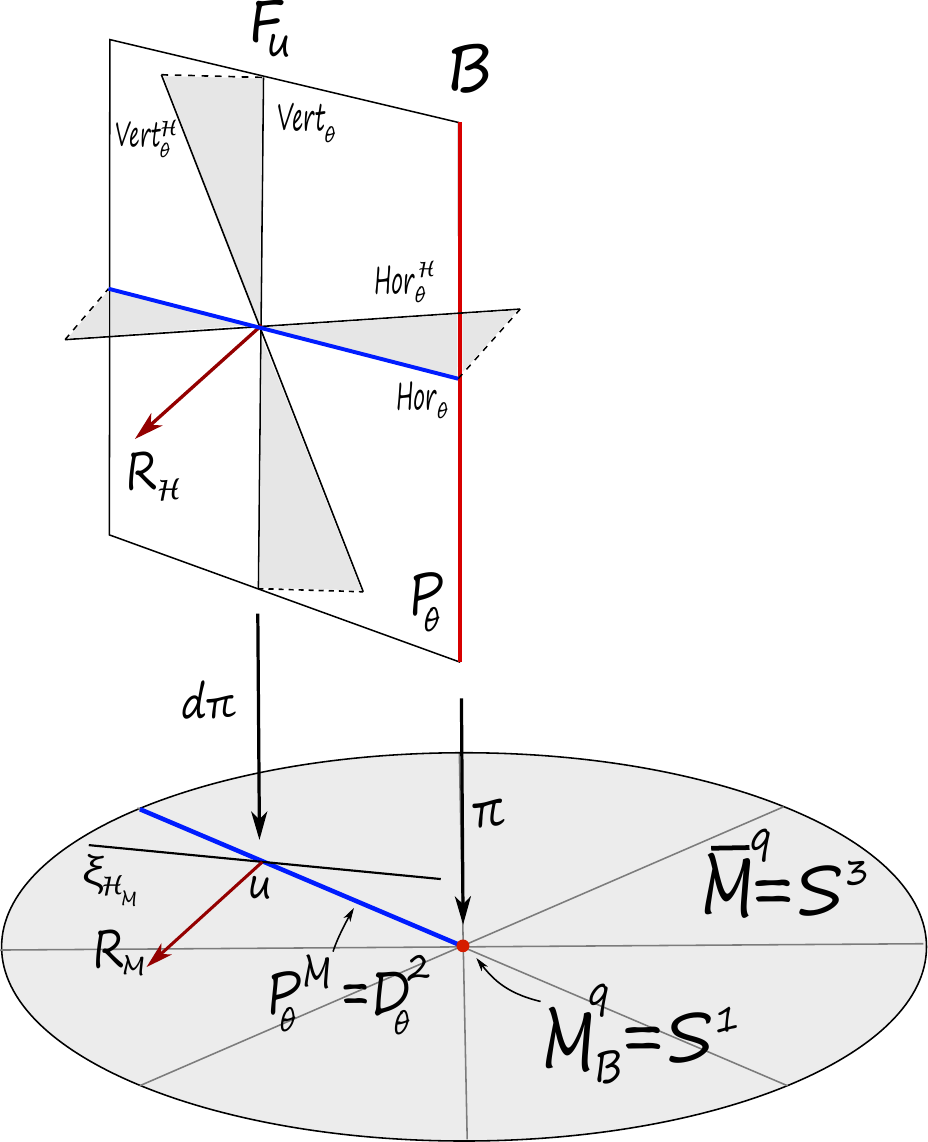}
    \caption{The complex splitting $TP_\theta=\mbox{Hor}_\theta \oplus \mbox{Vert}_\theta$, and the induced splitting $\xi_\mathcal{H}=\mbox{Hor}^\mathcal{H}_\theta \oplus \mbox{Vert}^\mathcal{H}_\theta$. The symplectic connection $\mbox{Hor}_\theta$ projects to $T_uP_\theta^\mathcal{M}$. Near the binding $B$, we have $\mbox{Hor}^\xi_\theta=\mbox{Hor}_\theta^\mathcal{H}$ and $\mbox{Vert}^\xi_\theta=\mbox{Vert}_\theta^\mathcal{H}$. Away from $B$, $\mbox{Hor}_\theta=\mbox{Hor}_\theta^\mathcal{H}$ and $\mbox{Vert}_\theta=\mbox{Vert}_\theta^\mathcal{H}$. This structure naturally induces a stable Hamiltonian structure deformation $\xi_{\mathcal{H}_\mathcal{M}}$ of $\xi_\mathcal{M}$ with Reeb vector field $R_\mathcal{M}$. The contact structure $\xi_\mathcal{M}$ is the \emph{average plane} in the cone $C_\xi=d\pi(\xi)$, and similarly, $\xi_{\mathcal{H}_\mathcal{M}}$ is the average plane in $C_\mathcal{H}=d\pi(\xi_\mathcal{H})$; see Lemma \ref{lemma:adapted} below.}
    \label{fig:horvert}
\end{figure}

\textbf{The complex normal bundle.} We can give the following explicit description of the complex normal bundle $N_u$ for each $u\in \overline{\mathcal{M}}$. Let $\mbox{Hor}_\theta$ denote the $\omega_\theta$-symplectic complement to $u$ inside $P_\theta$, which is a symplectic connection for $\pi_\theta$, i.e.\ $TP_\theta=\mbox{Hor}_\theta \oplus \mbox{Vert}_\theta$,  $\mbox{Vert}_\theta=\ker d\pi_\theta$, $\mbox{Hor}_\theta \cong T\mathbb{D}^2_\theta$ under $d\pi_\theta$; this splitting for $TP_\theta$ is also $J_\theta$-complex. Using the complex isomorphism  $(TP_\theta,J_\theta)\cong (\xi_\mathcal{H}\vert_{P_\theta}, J\vert_{\xi_\mathcal{H}\vert_{P_\theta}})$, the $J_\theta$-complex splitting $TP_\theta=\mbox{Hor}_\theta \oplus \mbox{Vert}_\theta$ then induces a $J$-complex splitting $\xi_\mathcal{H}\vert_{P_\theta}= \mbox{Hor}_\theta^\mathcal{H} \oplus \mbox{Vert}^\mathcal{H}_\theta$, and similarly a $J^\xi$-complex splitting $\xi =\mbox{Hor}_\theta^\mathcal{\xi} \oplus \mbox{Vert}^\mathcal{\xi}_\theta$. Under $d\pi: T(M\backslash L)\rightarrow T\overline{\mathcal{M}}^q$, $\mbox{Hor}_\theta$ projects to $TP_\theta^\mathcal{M}$.  See Figure \ref{fig:horvert}. If we denote by $\widehat N_\theta$ the $\omega$-symplectic complement of $\widehat P_\theta$, then, if $u\subset \widehat P_\theta$, we have $N_u=\mbox{Hor}_\theta\vert_u \oplus \widehat N_\theta\subset u^*T(\mathbb{R}\times M)$. Moreover, the tangent space $T_u\overline{\mathcal{M}}=\ker \mathbf{D}_u\subset W^{1,2}(N_u)$ is a $4$-dimensional subspace consisting of smooth sections of $N_u$, and we have the explicit description $\ker \mathbf{D}_u=\langle e_1,e_2 \rangle \oplus \langle \partial_t,n\rangle$, where $e_1,e_2=J_\theta(e_1)$ point-wise span $\mbox{Hor}_\theta\vert_u$ (or equivalently span $P_\theta^\mathcal{M}\vert_u$ when viewed as tangent to $\mathcal{M}^q$), and $n$ corresponds to the Hopf direction. In other words, vectors in $\mbox{Hor}_\theta$ correspond to nearby fibers of the Lefschetz fibration in the same page, $\partial_t$ corresponds to $\mathbb{R}$-translation, and $n$, to curves in nearby pages. We have an analogous description in the case when $u\subset \mathcal{M}_B$. 

\smallskip

\textbf{The shadowing open book supports the shadowing contact structure.} Since the open book $M=\mathbf{OB}(P,\phi)$ supports the contact structure $\xi$, away from the binding $B$ we have an isomorphism $\xi\vert_{P_\theta} \cong TP_\theta$, for each $\theta$. Therefore, the splitting $TP_\theta=\mbox{Vert}_\theta \oplus \mbox{Hor}_\theta$ induces a splitting $\xi\vert_{P_\theta} = \mbox{Vert}_\theta^\xi \oplus \mbox{Hor}_\theta^\xi$. Then $\xi_{\mathcal{M}}$ is given along the interior of the $\theta$-page $P_\theta^\mathcal{M}$ as the average 
$$\xi_{\mathcal{M}}\vert_{P_\theta^\mathcal{M}}=\int_{z\in F_u}d\pi^q(\mbox{Hor}^\xi_\theta(z))dz.$$
See Remark \ref{rk:integrals} below for details on how to order to interpret this expression. 

Since $\mbox{Hor}_\theta \cong T\mathbb{D}^2_\theta$, $\xi_{\mathcal{M}}\vert_{P_\theta^\mathcal{M}}$ is isomorphic to $TP_\theta^\mathcal{M}=d\pi^q(T\mathbb{D}^2_\theta)$. This, together with the fact that its binding $\mathcal{M}_B^q$ is a Reeb orbit (see the discussion of the holomorphic shadow below), means that the trivial open book supports the contact structure $\xi_\mathcal{M}$. Note that, for $u \in \mathcal{M}^q_B$, the contact structure $\xi_{\mathcal{M}}\vert_u$ is identified with $\xi_B^\perp\vert_u$. Here, $\xi_B^\perp$ is the symplectic normal bundle of $\xi_B$ inside $\xi\vert_B$ with respect to $d\alpha\vert_{\xi}$, which is symplectically trivial along $u$. One can then intepret the distribution $\mbox{Hor}^\xi$ across the binding as $\mbox{Hor}^\xi\vert_B=\xi_B^\perp$. 

\medskip

\textbf{The shadowing almost complex structure.} We have that $J^\xi$ induces a \emph{shadowing} almost complex structure $J^\xi_\mathcal{M}$ on $\ker \mathbf{D}_u$, as follows. Consider the Hilbert space $W^{1,2}(N_u)$ of $W^{1,2}$-sections of the normal bundle of $u$, endowed with the inner product 
$$
g_u(v,w)=\int_{z\in F_u} g_z(v(z),w(z))dz,
$$
where $g_z=d\alpha_z(\cdot,J^\xi \cdot)+\alpha_z\otimes \alpha_z+dt \otimes dt$. By varying $u$, we can view $g$ as a Riemannian metric on the Banach bundle $W^{1,2}(N_u)\hookrightarrow \mathcal{W} \rightarrow \mathcal{M}$. We denote by $g^q$ the metric induced by $g$ on the quotient Banach bundle $W^{1,2}(N_u\backslash \mathbb{R})\hookrightarrow \mathcal{W}^q=\mathcal{W}\backslash \mathbb{R} \rightarrow \mathcal{M}^q=\mathcal{M}\backslash \mathbb{R}$. The almost complex structure $J^\xi$ gives an endomorphsim
$$
J^\xi_u: W^{1,2}(N_u)\rightarrow W^{1,2}(N_u),
$$
$$
J^\xi_u(v(z))=J^\xi_z(v(z)),
$$
for $z\in F_u$, compatible with $g_u$. Then $\omega_u(v,w)=g_u(v,J^\xi_u(w))$ is a symplectic form on $W^{1,2}(N_u)$, which can be viewed as a fiber-wise symplectic form in the Banach bundle $\mathcal{W}$. Let $$P_u: W^{1,2}(N_u)\rightarrow \ker \mathbf{D}_u$$ denote the $L^2$-orthogonal projection with respect to $g_u$. Then $J^\xi_\mathcal{M}$ is defined via $$(J^\xi_\mathcal{M})_u(v(u))(z)=P_u(J_u^\xi(v(z))),$$ for $z\in F_u$, and $v\in \ker \mathbf{D}_u$. It is compatible with $\omega_\mathcal{M}^\varphi$. The splitting $\ker \mathbf{D}_u=\xi_\mathcal{M}\oplus \langle \partial_t,R_\mathcal{M}\rangle$ is $J_\mathcal{M}^\xi$-complex, where $R_\mathcal{M}$ denotes the Reeb vector field of $\alpha_\mathcal{M}$, which follows from the fact that $P_u(R_\alpha\vert_u)$ is a multiple of $R_\mathcal{M}$ (see Equation (\ref{eq:L2rep}) below).

\medskip

\textbf{The shadowing stable Hamiltonian structure.} We have a \emph{shadowing} stable Hamiltonian structure $\mathcal{H}_\mathcal{M}=(\lambda_\mathcal{M},d\alpha_\mathcal{M})$ on $\overline{\mathcal{M}}^q$ whose kernel $\xi_{\mathcal{M}_H}$ deforms $\xi_\mathcal{M}$, exactly as for $\mathcal{H}$ and $\xi$, i.e.\ $\xi_{\mathcal{H}_\mathcal{M}}$ is tangent to $TP_\theta^\mathcal{M}$ away from $\mathcal{M}_B^q$, and the Reeb vector field of $\mathcal{H}_\mathcal{M}$ is $R_\mathcal{M}$, the Reeb vector field of $\alpha_\mathcal{M}$. The splitting $\ker \mathbf{D}_u=\mathcal{H}_\mathcal{M}\oplus \langle \partial_t,R_\mathcal{M}\rangle$ is $J_\mathcal{M}$-complex, where $J_\mathcal{M}$ is induced from $J$ analogously as above, i.e.\ via $(J_\mathcal{M})_u(v(z))=P^\mathcal{H}_u(J_z(v(z)))$, where $P^\mathcal{H}_u: W^{1,2}(N_u)\rightarrow \ker \mathbf{D}_u$ is the orthogonal projection with respect to 
$$
g^\mathcal{H}_u(v,w)=\int_{z\in F_u} g^\mathcal{H}_z(v(z),w(z))dz,
$$
with $g^\mathcal{H}_z=d\alpha_z(\cdot,J \cdot)+\lambda_z\otimes \lambda_z+dt \otimes dt$.

\medskip

\textbf{The (pseudo-)holomorphic shadow.} We call the Reeb flow of $\alpha_\mathcal{M}$ on $\overline{\mathcal{M}}^q$, the \emph{holomorphic shadow} of the Reeb flow of $\alpha$ on $M$, generated by $R_\mathcal{M}$. In fact, $R_{\mathcal{M}}\vert_u$ is a positive reparametrization of the $L^2$-orthogonal projection of $R_\alpha\vert_u$ (a vector field in $M$ normal to $u$, i.e.\ in $W^{1,2}(N_u)$) to the kernel of $\mathbf{D}_u^N$, i.e.\ the tangent space to the moduli space. This can be proved as follows.

By construction, the Reeb vector field $R_\mathcal{M}\in T\overline{\mathcal{M}}^q$ of $\alpha_\mathcal{M}$ is defined by the equations
$$
\mathbf{D}_uR_\mathcal{M}=0,
$$
\begin{equation}\label{eq:Reebeqs}
1=(\alpha_\mathcal{M})_u(R_\mathcal{M}(u))=\int_{z\in F_u}\alpha_z(R_\mathcal{M}(z))dz,
\end{equation}
$$
0=(d\alpha_\mathcal{M})_u(R_\mathcal{M}(u),\cdot)=\int_{z\in F_u}d\alpha_z(R_\mathcal{M}(z),\cdot)dz.
$$
We claim that
\begin{equation}\label{eq:L2rep}
R_\mathcal{M}(u)=\frac{P_u(R_\alpha\vert_u)}{(\alpha_\mathcal{M})_u(P_u(R_\alpha\vert_u))}.
\end{equation}
This can be checked as follows. Note that the metric $g^q_\mathcal{M}=g^q\vert_{T\mathcal{M}^q}$, where we view $T\overline{\mathcal{M}}^q\subset \mathcal{W}^q$ as a rank-$3$ subbundle, can in fact be expressed as 
$$
g^q_\mathcal{M}=d\alpha_\mathcal{M}(\cdot,J^\xi_\mathcal{M}\cdot)+\Vert R_\mathcal{M}\Vert^2\alpha_\mathcal{M}\otimes \alpha_\mathcal{M},
$$
where 
$$
\Vert R_\mathcal{M}(u) \Vert^2=(g^q_\mathcal{M})_u(R_\mathcal{M}(u),R_\mathcal{M}(u))=\int_{z\in F_u} \left[\alpha_z(R_\mathcal{M}(z))\right]^2dz.
$$
We then consider a basis $e_1,e_2=J^\xi_\mathcal{M}(e_1) \in \xi_{\mathcal{M}}, e_3=\frac{1}{\Vert R_\mathcal{M}\Vert}R_\mathcal{M}$ of $T_u\mathcal{M}^q$, which is orthonormal with respect to $g^q_\mathcal{M}$. Then $$
P_u(R_\alpha\vert_u)=\sum_{i=1}^3g_u(R_\alpha\vert_u,e_i(u))e_i(u).
$$

Moreover, 
$$
g_u(R_\alpha\vert_u,e_i(u))=\int_{z\in F_u} g_z(R_\alpha(z),e_i(z))dz=\int_{z\in F_u} \alpha_z(e_i(z))dz=(\alpha_\mathcal{M})_u(e_i(u)),
$$
and therefore
$$P_u(R_\alpha\vert_u)=\frac{R_\mathcal{M}(u)}{\Vert R_\mathcal{M}(u)\Vert^2},
$$
with $(\alpha_\mathcal{M})_u(P_u(R_\alpha\vert_u))=\frac{1}{\Vert R_\mathcal{M}(u)\Vert^2}$. This proves the claim.

\smallskip

\textbf{The shadow as a flow on $M$ preserving the foliation.} To gain some more insight on the holomorphic shadow, we need better understanding of the tangent space to $\overline{\mathcal{M}}^q$. If $u$ lies in a hypersurface $H \in \mathcal{L}$, there is a natural splitting of the normal operator $\mathbf{D}_u^N=\mathbf{D}_u^{T_H} \oplus \mathbf{D}_u^{N_H}$ into tangent and normal components with respect to $H$ (cf.\ \cite[sec.\ 4.7]{Mo}). Here, $N_H$ is the normal bundle to $H$, which is a choice of holomorphic complement to $H$ which near infinity is $\mathbb{R}$-invariant and coincides with $\xi_B^\perp \subset \xi\vert_B$. If $H$ is not cylindrical and projects to the $\theta$-page $P_\theta$, then $\ker \mathbf{D}_u^{N_H}$ is $2$-dimensional, spanned by the $\mathbb{R}$-direction $\partial_t$ and a vector field $n_M\vert_u$ which takes values in $N_H\vert_u$; while $\ker \mathbf{D}_u^{T_H}$ is identified under $d\pi_\theta\vert_u$ with $T_{\pi_\theta(u)}\mathbb{D}^2_\theta$, parametrizing the fibers of the Lefschetz fibration $\pi_\theta$ near to $u$. If $H=\mathcal{L}_B$ is cylindrical over $B$, then $\ker \mathbf{D}_u^{N_H}=0$ while $\ker \mathbf{D}_u^{T_H}$ is spanned by $\partial_t$ and a vector field $n_B\vert_u$ which takes values in $N_u\cap TB$. We may choose the vector fields $n_B$ and $n_M$ so that they glue smoothly together to a vector field $n$, in such a way that, when we identify $\mathcal{M}^q\cong S^3$, the flow of $n$ is precisely the Hopf flow on $S^3$. The projection of $R_\alpha\vert_u$ to $\ker \mathbf{D}_u^N$ can then be written as
$$
R_{\mathcal{M}}\vert_u=\left\{\begin{array}{cc}
  F(u)n\vert_u,   &  \mbox{ if } u \in \mathcal{M}_B^q \\
    F(u)n\vert_u+R_{\mathcal{M}}^T\vert_u & \mbox{ if } u \in \overline{\mathcal{M}}^q\backslash \mathcal{M}_B^q
\end{array}\right.,
$$
where $F: \overline{\mathcal{M}}^q\rightarrow \mathbb{R}^+$ is a smooth positive function, and $R_{\mathcal{M}}^T\vert_u \in \ker \mathbf{D}^{T_H}_u\cong T_{\pi_\theta(u)} \mathbb{D}^2$, where $H \in \mathcal{H}$ projects to $P_\theta$ and contains $u$. This fully completes the proof that the trivial open book in $\overline{\mathcal{M}}^q=\mathbf{OB}(\mathbb{D}^2,\mathds{1})$ is adapted to the Reeb dynamics of $\alpha_\mathcal{M}$, its binding $\mathcal{M}^q_B$ being a Reeb orbit. This finishes the proof of Theorem \ref{thm:contactstructure}.

\smallskip

\textbf{The shadowing cone.} We define the \emph{shadowing cone} by $C_\alpha:=d\pi(\ker d\alpha)\subset T\overline{\mathcal{M}}^q$, where $\pi: M\backslash L \rightarrow \overline{\mathcal{M}}^q$ is the quotient map to the leaf space. We also define the cones $C_\xi:=d\pi(\xi)$, and $C_\mathcal{H}:=d\pi(\xi_\mathcal{H})$.

\begin{lemma}\label{lemma:adapted}
The shadowing cone $C_\alpha$ is strongly adapted to the trivial open book $(\mathcal{M}_B^q,\theta_\mathcal{M})$ in $\overline{\mathcal{M}}^q$. Namely, we have the following:
\begin{enumerate}
    \item The contact form $\alpha_\mathcal{M}$ is a section of $C_\alpha$;
    \item The $1$-form $d\theta_\mathcal{M}$ is a section of $C_\alpha$;
    \item The $1$-dimensional cone $\ker d\alpha_\mathcal{M}$ is a subcone of $C_\alpha$, and in fact is the average direction in $C_\alpha$ (i.e.\ $C_\alpha$ is centered at $R_\mathcal{M}$). We write
    \begin{equation}\label{averagedirection}
           \ker d\alpha_\mathcal{M}=\int_{z\in F_u} (d\pi\ker d\alpha_z)dz.
    \end{equation}
    In particular $R_\mathcal{M}\in C_\alpha$;
    \item $C_\alpha\vert_{\mathcal{M}^q_B}= T\mathcal{M}^q_B$;
    \item Similarly, the contact structure $\xi_\mathcal{M}$ is the average plane in $C_\xi$, and $\xi_{\mathcal{H}_\mathcal{M}}$, the average plane in $C_\mathcal{H}$, i.e.\
\begin{equation}\label{averageplanes}
(\xi_\mathcal{M})_u=\int_{z\in F_u} d\pi(\xi_z)dz,\;(\xi_\mathcal{H_\mathcal{M}})_u=\int_{z\in F_u} d\pi((\xi_\mathcal{H})_z)dz.
\end{equation}
\end{enumerate}

\begin{remark}\label{rk:integrals}
Formula (\ref{averagedirection}) can be interpreted as follows. The intersection of $C_\alpha$ with the $g^q_\mathcal{M}$-unit sphere $S_\mathcal{M}\subset T\overline{\mathcal{M}}^q$, which parametrizes the directions in $C_\alpha$, is, away from $\mathcal{M}^q_B$, a closed subset of the upper hemisphere $S_\mathcal{M} \cap \{d\theta_\mathcal{M}> 0\}$ of $S_\mathcal{M}$. This subset can therefore be globally parametrized by two angles $(\phi_1,\phi_2)$ (longitude and latitude); the direction $\ker d\alpha_\mathcal{M}$ then corresponds to $$(\phi_1(R_\mathcal{M}),\phi_2(R_\mathcal{M}))=\frac{R_\mathcal{M}}{\Vert R_\mathcal{M} \Vert}=\int_{z\in F_u}(\phi_1(z),\phi_2(z))dz.$$ Along $\mathcal{M}_B^q$, Equation (\ref{averagedirection}) also holds, where the average direction is the only direction, i.e.\ $$\ker d\alpha_\mathcal{M}\vert_{\mathcal{M}_B^q}=T^+\mathcal{M}^q_B=C_\alpha\vert_{\mathcal{M}_B^q},$$ where $T^+\mathcal{M}^q_B$ are the vectors in $T\mathcal{M}^q_B$ which point in the positive direction according to the orientation of $\mathcal{M}^q_B=S^1$. The expressions in (\ref{averageplanes}) can be understood in a similar fashion. Namely, the metric $g^q_\mathcal{M}$ induces a diffeomorphism between the Grassmannian of co-oriented $2$-planes in $T\overline{\mathcal{M}}^q$ and $S_\mathcal{M}$, by mapping such a $2$-plane to the unique unit vector which is positively orthogonal to the plane; the above considerations then apply.
\end{remark}

\end{lemma}
\begin{proof}[Proof of Lemma \ref{lemma:adapted}]
For the first part, we need to show that $(\alpha_\mathcal{M})_u(d_z\pi(R_\alpha(z)))$ is positive for $z\in F_u$, $u=\pi(z)$. If $e_1,e_2=J^\xi_\mathcal{M}(e_1) \in \xi_{\mathcal{M}}, e_3=R_\mathcal{M}$ is the $g^q_\mathcal{M}$-orthogonal basis for $T_u\mathcal{M}^q$ considered above, we have
$$
d_z\pi(R_\alpha(z))=\sum_{i=1}^3 g_z(R_\alpha(z),e_i(z))e_i(u)=\sum_{i=1}^3 \alpha_z(e_i(z))e_i(u),
$$
and so
$$
(\alpha_\mathcal{M})_u(d_z\pi(R_\alpha(z)))=\alpha_z(R_\mathcal{M}(z)),
$$
which is strictly positive. This proves the first claim. The second claim follows immediately from the fact that $R_\alpha$ is positively transverse to $P_\theta$ for every $\theta$, away from $B$. The third and fifth follows by construction. The fourth follows from the fact that $B$ is invariant under $R_\alpha$. This proves the lemma. 
\end{proof}

\begin{example}[Sasakian case=Kepler problem]\label{rk:Sasakian} Note that while $R_\alpha$ satisfies the second and third of the Equations (\ref{eq:Reebeqs}), it might not satisfy the first one; this is true if e.g.\ the flow of $R_\alpha$ preserves $J$, i.e.\ it is holomorphic, which holds if the strict contact manifold $(M,\alpha)$ is $K$-contact or Sasakian. For instance, the unit cotangent bundle of a Riemannian manifold $(X,g)$ carries a Sasakian structure with contact form induced from the standard Liouville form, and metric $\tilde g$ induced from $g$, via its Levi-Civita connection. In particular, we have the case of $(M,\alpha)=(S^*S^3,\alpha_g)$ where $g$ is the standard round metric and $\alpha_g$ the standard Liouville form, whose flow is the round geodesic flow on $S^3$. This is an IP $5$-fold, and $(S^*S^3,\alpha_g)$ is Sasakian, its Reeb flow preserving the metric $\tilde g$ and the almost complex structure on $\ker \alpha_g$ induced from the restriction of the integrable complex structure on the Stein manifold $T^*S^3$. For this case, we have $R_\mathcal{M}=R_\alpha$, inducing the Hopf flow on $S^3$. This dynamical system corresponds to the spatial Kepler problem after Moser regularization. The holomorphic shadow for the spatial Kepler problem is then the Hopf flow.
\end{example}

\begin{example}[Rotating Kepler problem]\label{ex:rotKepler}
In \cite[Appendix A]{MvK}, the author and Otto van Koert constructed a symplectic Lefschetz fibration on the page $\mathbb{D}^*S^2$ of an open book, whose fibers are annuli, and which are invariant under the return map for the rotating Kepler problem. Its shadow is then also the Hopf flow, and the shadowing cone is its positive span. Heuristically speaking, when the mass ratio $\mu$ is then perturbed to be small, the expectation is that the shadowing cone then has non-empty interior, and its size provides some sort of ``measure'' of non-integrability of the problem. 
\end{example}

\textbf{A non-trivial example: the Katok example in dimension $5$.} Now we look at the $5$-dimensional instances of the well-known examples by Katok \cite{K73}, of Finsler metrics on $S^3$ with only finitely many simple closed geodesics. It turns out that the holomorphic shadow construction can be understood explicitly, and corresponds to the Reeb dynamics on the boundary of suitable irrational ellipsoids. 

We follow the discussion in Appendix A in \cite{MvK2}. Namely, we consider the Brieskorn manifold
$$
\Sigma^{5} := \left\{ (z_0,z_1,z_2,z_3) \in \C^{4} ~\Bigg|~\sum_j z_j^2 = 0 \right\}
\cap 
S^7,
$$
equipped with the contact form $\alpha =\frac{i}{2}\sum_j z_j d\bar z_j -\bar z_j dz_j$.
This space is contactomorphic to $S^*S^3$ with its canonical contact structure. We consider the unitary change of coordinates:
$$
w_0=z_0,
w_1=z_1,
w_{2}=\frac{\sqrt 2}{2} (z_{2}+i z_{3}),
w_{3}=\frac{i \sqrt 2}{2} (z_{2}-i z_{3}).
$$

For $\epsilon>0$ small and irrational, define the function
$$
H_\epsilon(w)=\Vert w \Vert^2+\epsilon( 
|w_{2}|^2-|w_3|^2
),
$$
and perturb the contact form as
$\alpha_\epsilon = H_\epsilon^{-1} \cdot \alpha.$
The Reeb flow is seen to be 
$$
(w_0,w_1,w_2,w_3)\longmapsto
(e^{2\pi it}w_0,e^{2\pi it}w_1,e^{2\pi it(1+\epsilon)}w_2,e^{2\pi it(1-\epsilon)}w_3).
$$
This flow has only $4$ periodic orbits by irrationality of $\epsilon$. We have a supporting open book for the contact form $\alpha_\epsilon$ given by
$$
\pi_0: \Sigma^{5} \longrightarrow \C,
$$
$$
(w_0,w_1,w_2,w_3) \longmapsto w_0.
$$
The zero set of $\pi_0$ defines the binding $\mathbb{R}P^3$, and the pages are the sets of the form $P_\theta=\{\arg \pi_0=\theta\}$, $\theta \in S^1$, which are all copies of $\mathbb{D}^*S^{2}$. As in \cite{MvK2}, only two of the four orbits correspond to interior fixed points of the resulting return map, which actually has no boundary fixed points; the other two orbits lie in the binding, but have periods which differ from the return time. Similarly, define the maps
$$
\pi_i: \Sigma^5\rightarrow \mathbb{C},
$$
$$
(w_0,w_1,w_2,w_3) \longmapsto w_i,
$$
for $i=1,2,3$. The Reeb flow then preserves the set of fibers of $\pi_i$. For $i=1$, $\pi_1$ is the standard Lefschetz fibration on each page $P_\theta=\mathbf{LF}(\mathbb{D}^*S^1,\tau_P^2)$, where $\tau_P$ is the Dehn twist. We will refer to the flow induced on the fiber-space of $\pi_i$ as the $i$-th holomorphic shadow.

Recall that the ellipsoid $E(a,b)$ is defined as
$$
E(a,b)=\left\{(u,v)\in \mathbb{C}^2: \frac{\pi|u|^2}{a}+\frac{\pi|u|^2}{b}\leq 1\right\},
$$
which is a star-shaped domain in $\mathbb{C}^2$ and so its boundary $S(a,b):=\partial E(a,b)\cong S^3$, which we call a \emph{spheroid} to avoid confusion (although this is also usually called an ellipsoid in the literature), inherits a contact form whose Reeb flow is $\phi_t(u,v)=(e^{2\pi i a t}u,e^{2 \pi i bt}v)$. One may view this as a flow on $S^3=S(1,1)$, which is adapted to the trivial open book given by $(u,v)\mapsto \frac{v}{|v|}$, and the return map on the disk-like pages is simply the rotation by angle $2\pi\frac{a}{b}$. We say that the ellipsoid/spheroid is irrational if $a$ and $b$ are rationally independent. In this case, the first return map at each page has the origin as the unique interior fixed point, and no boundary fixed points. 

From the above discussion, we easily conclude:

\begin{proposition}\label{prop:katok}
The $i$-th holomorphic shadow of the $5$-dimensional Katok example is the Reeb dynamics of:
\begin{itemize}
    \item $S(1,1)$, i.e.\ the Hopf flow, for $i=1$;
    \item The irrational spheroid $S(1,1+\epsilon)$, for $i=2$;
    \item The irrational spheroid $S(1,1-\epsilon)$, for $i=3$. 
\end{itemize}
\end{proposition}

\medskip
One can similarly choose the open book to be a $\pi_i$, and the rest of the $\pi_j$ to be the Lefschetz fibrations on a fixed page of $\pi_i$. The dynamics on the irrational spheroid $S(1+\epsilon,1-\epsilon)$ is then the $3rd$-holomorphic shadow with respect to $\pi_2$. We remark that, rather interestingly, this is precisely the one which doubly covers the $3$-dimensional instance of Katok's construction on $\mathbb{R}P^3$ \cite[Thm.\ 2.3, rk.\ 2.17]{AGZ} (and the latter we see in the binding of $\pi_0$).  

\medskip

\textbf{Proof of Lifting theorem.} We prove surjectivity of $\mathbf{HS}$.

\begin{proof}[Proof of Thm.\ \ref{thm:lifting}] The proof is straightforward from the construction of $\mathbf{HS}$, but it requires some care as to what is fixed and what is a choice, making it rather involved. For convenience of the reader, we provide careful details. 

Denote by $\alpha_{0}$ the contact form on $S^3$ whose Reeb dynamics is adapted to the concrete open book $\pi_0:S^3\backslash L\rightarrow S^1$, where $L \cong S^1$, having disk-like $\theta$-page $P_\theta^0$. Let $\xi_0=\ker \alpha_0$ be the induced contact distribution.  

We fix the following choices:
\begin{itemize}
    \item A contact form $\alpha_B \in \mathbf{Reeb}(F,\phi_F)$ (e.g.\ via the Thurston-Wilkelnkemper construction of a Giroux form), adapted to some choice of concrete open book on $B$;
    \item An exact symplectic form $\omega_\theta=d\lambda_\theta$ on each page $P_\theta$ in such a way that $(P_\theta,\omega_\theta)$ is an ideal Liouville filling of $(B,\alpha_B)$ for every $\theta$;
    \item An $\omega_\theta$-compatible almost complex structure $J_\theta$ on $P_\theta$, which is generic as a $1$-parameter family.
    \item An almost complex structure $J_B$ on $\mathbb{R}\times M$ which makes the concrete open book on $B$ holomorphic.
\end{itemize}
This data can be used to construct a Giroux form $\alpha \in \mathbf{Reeb}(\mathbf{LF}(F,\phi_F),\phi)$, satisfying $\alpha_B=\alpha\vert_B$ and $\omega_\theta=d\alpha\vert_{P_\theta}$ for each $\theta$, which is \emph{auxiliary}. This gives an auxiliary concrete contact distribution $\xi=\ker \alpha$ on $M$, supported by the IP open book. One can also construct an almost complex structure $J$, such that $J\vert_{\xi\vert_{P_\theta}}$ corresponds to $J_\theta$ under the isomorphism $\xi\vert_{P_\theta}\cong TP_\theta$, and such that $J\vert_B=J_B$; we use this $J$ to construct the holomorphic foliation $\overline{\mathcal{M}}^q$ associated to $\alpha$.

Having made these choices, away from $B$, we have a symplectic splitting
$$
(TP_\theta,\omega_\theta)=(\mbox{Vert}_\theta,\omega_\theta^v)\oplus(\mbox{Hor}_\theta,\omega_\theta^h):=(\ker d\pi_\theta,\omega_\theta\vert_{\ker d\pi_\theta})\oplus (\mbox{Vert}_\theta^{\perp_{\omega_\theta}},\omega_\theta\vert_{\mbox{Vert}_\theta^{\perp_{\omega_\theta}}}),
$$
together with a symplectic isomorphism $(\xi\vert_{P_\theta},d\alpha\vert_{\xi\vert_{P_\theta}})\cong (TP_\theta,\omega_\theta)$. This induces a symplectic splitting
$$
(\xi\vert_{P_\theta},d\alpha\vert_{\xi\vert_{P_\theta}})=(\mbox{Vert}^{\xi}_\theta,\omega_\theta^{v,\xi})\oplus(\mbox{Hor}^{\xi}_\theta,\omega_\theta^{h,\xi}),
$$
where $\omega_{\theta}^{v,\xi}=d\alpha\vert_{\mbox{Vert}^{\xi}_\theta},$ $\omega_\theta^{h,\xi}=d\alpha\vert_{\mbox{Hor}^{\xi}_\theta}$. The symplectic distribution $(\mbox{Hor}^{\xi},\omega^{h,\xi})$ extends to $B$ as $(\xi_B^{\perp_M},d\alpha\vert_{\xi_B^{\perp_{M}}})$, where the symplectic complement is with respect to $d\alpha$.

The holomorphic foliation $\overline{\mathcal{M}}^q$ and concrete open book on $M$ provide a diffeomorphism $\mathcal{F}: S^3 \rightarrow \overline{\mathcal{M}}^q$, so that we have an induced concrete trivial open book $\pi_\mathcal{M}:\overline{\mathcal{M}}^q\backslash \mathcal{M}^q_B\rightarrow S^1$ with $\theta$-page $P_\theta^\mathcal{M}$ together with a contact form $\alpha_{\mathcal{M}}\in \Omega^1(\overline{\mathcal{M}}^q)$ and a concrete contact distribution $\xi_{\mathcal{M}}=\ker \alpha_{\mathcal{M}}$ supported by the open book. Via the pullback connection  $\mbox{Hor}_\theta^\mathcal{M}=(ev^q)^*(\mbox{Hor}_\theta)\subset T\overline{\mathcal{M}}^q_*\vert_{\pi_*^{-1}(P_\theta^{\mathcal{M}})}$ for the forgetful map $\pi_*$, away from $B$ we have an induced identification $TP_\theta^\mathcal{M}\vert_u=\mbox{Hor}_\theta\vert_u$ as a two-plane distribution on $M$ along $u$. We also have a two-plane distribution $\mbox{Hor}_\theta^{\xi_\mathcal{M}}\subset TM\vert_{P_\theta}$, which is identified with $\xi_\mathcal{M}^\theta:=\xi_\mathcal{M}\vert_{P_\theta^\mathcal{M}}$, as well as, away from $B$, an induced linear isomorphism $\mbox{Hor}_\theta\cong \mbox{Hor}_\theta^{\xi_\mathcal{M}}$ covering the identity, coming from the analogous isomorphism $TP_\theta^\mathcal{M}\cong \xi_\mathcal{M}^\theta$ away from $\mathcal{M}_B^q$. This induces a symplectic form $\omega_\theta^{h,\xi_{\mathcal{M}}}$ on $\mbox{Hor}_\theta^{\xi_\mathcal{M}}$, coming from $\omega_\theta^{h}$. Moreover, we can view $d\alpha_{\mathcal{M}}$ as a non-degenerate $2$-form on $\mbox{Hor}^{\xi_\mathcal{M}}$. For $u \in \mathcal{M}_B^q$, $\xi_{\mathcal{M}}\vert_u$ gets identified with some $2$-plane distribution  $\xi_B^{\perp_{\mathcal{M}}}\subset TM\vert_B$ transverse to $TB$ which also comes with the non-degenerate $2$-form $d\alpha_\mathcal{M}$. Note that $\mbox{Hor}_\theta^{\xi_\mathcal{M}}$ might not necessarily be a subbundle of the concrete contact distribution $\xi=\ker \alpha$, i.e.\ it does not necessarily agree with $\mbox{Hor}_\theta^{\xi}$. To remedy this, we may take an isotopy $\{\alpha^t\}_{t \in [0,1]}$ of Giroux forms, with $\alpha^0=\alpha$, satisfying $d\alpha^t\vert_{P_\theta}=\omega_\theta$ and $\alpha^t\vert_{B}=\alpha_B$ independently of $t$, together with the corresponding splitting
$$
(\xi^t\vert_{P_\theta},d\alpha^t\vert_{\xi^t\vert_{P_\theta}})=(\mbox{Vert}^{\xi^t}_\theta,\omega_\theta^{v,\xi^t})\oplus(\mbox{Hor}^{\xi^t}_\theta,\omega_\theta^{h,\xi^t}),
$$
so that $\mbox{Hor}^{\xi^1}=\mbox{Hor}^{\xi_\mathcal{M}}$ (and in particular $\xi_B^{\perp_M}=\xi_B^{\perp_\mathcal{M}}=:\xi_B^\perp$ along $B$). This induces an accompanying homotopy of almost complex structures $\{J_t\}$, so that $J_t\vert_B=J_B$, and $J_t$ is independent of $t$ under the isomorphism $\xi^t\vert_{P_\theta}\cong TP_\theta$ (i.e.\ it corresponds to $J_\theta$). The corresponding foliations $\overline{\mathcal{M}}^q_t$ are then $t$-independent.  

Up to this modification, we have $d\alpha_{\mathcal{M}}=f \omega^{h,\xi}$ for a unique positive smooth function $f:M\backslash L \rightarrow \mathbb{R}^+$. Since $\mathcal{M}^q_B$ is invariant under the flow of $R_\mathcal{M}$, we have $df\vert_{\xi_B^\perp}=0$. We then consider the contact form $\alpha^f:=f\alpha$, and its Reeb vector field $R_\alpha^f$. By construction, we have $\ker \alpha^f=\xi$, $d\alpha^f\vert_{\mbox{Hor}^{\xi}_\theta}=d\alpha_{\mathcal{M}}\vert_{\mbox{Hor}^{\xi}_\theta}$, and $R_\alpha^f\vert_B\in TB$ so that $B$ is invariant under the flow of $R_\alpha^f$.

Now, the Reeb vector field $R_{\mathcal{M}}$ of $\alpha_{\mathcal{M}}$ can be viewed as the unique vector field on $M$ characterized by the equations:
\begin{enumerate}
    \item $\mathbf{D}_u^N R_{\mathcal{M}}\vert_u=0$ for each $u \in \overline{\mathcal{M}}$;
    \item $d\alpha_{\mathcal{M}}(R_\mathcal{M},\cdot)=0$;
    \item $\alpha_{\mathcal{M}}(R_\mathcal{M})=1$.
\end{enumerate}
While $R_{\mathcal{M}}$ is positively transverse to $TP_\theta,$ as well as to $\mbox{Hor}_\theta^{\xi}$, a priori it might fail to be positively transverse to $\mbox{Vert}_\theta^{\xi}$. Since we still have freedom in the vertical directions, we can then further modify the Giroux form $\alpha$ via an isotopy as above, without changing the horizontal distribution and relative $B$, and without changing the foliation $\overline{\mathcal{M}}^q$, so that $\eta \in T\overline{\mathcal{M}}^q$ satisfies $\alpha_{\mathcal{M}}(\eta)>0$ if and only if $\alpha(\eta)>0$ pointwise, when viewed as a vector field in $M$. If $P_u: W^{1,2}(N_u)\rightarrow \ker \mathbf{D}_u^N$ denotes the orthogonal projection with respect to the metric induced by $J$ and $d\alpha$, the projection $P_u(R^f_\alpha\vert_u)$ of $R^f_\alpha\vert_u$ satisfies the first two equations above, and $g\vert_u:=\alpha_\mathcal{M}(P_u(R^f_\alpha\vert_u))>0$ defines a smooth positive function $g:M\backslash L \rightarrow \mathbb{R}^+$. Then $R_\mathcal{M}=g^{-1}P_u(R^f_\alpha\vert_u)$, so that $R_{\mathcal{M}}$ is the holomorphic shadow of $R^f_\alpha$. Note that $R^f_\alpha$ is positively transverse to the pages: if it were tangent to a page at a point, since pages are invariant under $J$ and symplectic for $d\alpha$, $R_{\mathcal{M}}$ would also be tangent to the page at that point, which is a contradiction. Then $\alpha^f \in \mbox{Reeb}(P,\phi)$ satisfies $\mathbf{HS}(\alpha^f,J)=\alpha_{\mathcal{M}}$.
\end{proof}

\section{Transverse paths and symplectic tomographies}\label{sec:tomographies}

\begin{figure}
    \centering
    \includegraphics[]{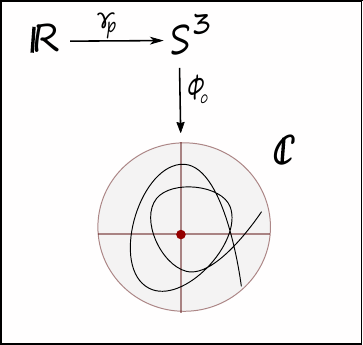}
    \caption{The qualitative image of a transverse path $\gamma_p$ under the defining map $\Phi_0:S^3\rightarrow \mathbb{C}$ for the trivial open book $\pi_0=\Phi_0/\vert \Phi_0\vert$.}
    \label{fig:transversepath}
\end{figure}

Alternatively to the construction of the holomorphic shadow, we may keep track of which holomorphic curves are intersected by each Reeb orbit in $M$, without changing the original dynamics. Namely, for $p \in M\backslash L$, we may consider the path
$$
\gamma_p(t)=\pi^q(\phi_t^M(p))\in \overline{\mathcal{M}}^q,
$$
where $\phi_t^M$ is the flow of $R_\alpha$. If $p \in B\backslash L$, this is a parametrization of $\mathcal{M}^q_B \cong S^1$; if $p \in M\backslash B$, this is a path in $S^3$ which is positively transverse to each disk-like page $P_\theta^\mathcal{M}$ as well as to the contact structure $\xi_\mathcal{M}$, and is by construction an orbit of the shadowing cone $C_\alpha$. Note that different choices of $p$ might induce paths which intersect each other (corresponding to their orbits intersecting the same holomorphic curve), and even self-intersect (corresponding to an orbit intersecting the same holomorphic curve multiple times), so these paths are not orbits of an autonomous flow. See Figure \ref{fig:transversepath}. We will refer to the collection
$$
\mathbf{TS}(\alpha,J)=\{\gamma_p: p \in M\backslash L\}
$$
as the \emph{transverse shadow} of the Reeb flow of $\alpha$ on $M$, with respect to $J$, which is by definition the collection of those orbits of $C_\alpha$ coming from orbits of $\alpha$ on $M$.

One may further choose to ``package'' these transverse paths in different ways, e.g.\ by considering those paths induced by points on a section of the Lefschetz fibration at a given page, as follows. 

Fix the $0$-page $P_0$ in $M$, and consider a two disk $D \cong \mathbb{D}^2$ satisfying:
\begin{itemize}
    \item $D\subset P_{0}$;
    \item $\partial D \subset \partial P_{0}=B$ is a loop which is disjoint from the binding $L$ of the concrete open book in $B$, and transverse to the interior of each of its pages and to the contact structure $\xi_B$;
    \item $D_0$ is a symplectic section of the Lefschetz fibration $\pi_0$, i.e.\ $D$ intersects each fiber of $\pi_0$ precisely once, and hence $D=im(s)$ for $s: \mathbb{D}_0^2\rightarrow P_0$ satisfying $\pi_0 \circ s=id$.
\end{itemize}

We refer to such a disk $D$ as a \emph{(horizontal) symplectic tomography} for the Reeb dynamics on $M$. Note that, if $\partial D$ is a Reeb orbit of $\alpha_B$ which is linked with $L$, then $f(\partial D)=\partial D$ is invariant under the return map. This is closely related to the counts of holomorphic sections with Lagrangian boundary condition of a given Lefschetz fibration, as considered e.g.\ by Seidel in \cite{Sei08}, with the difference that we consider the more flexible class of symplectic ones.

\begin{figure}
    \centering
    \includegraphics[width=0.6 \linewidth]{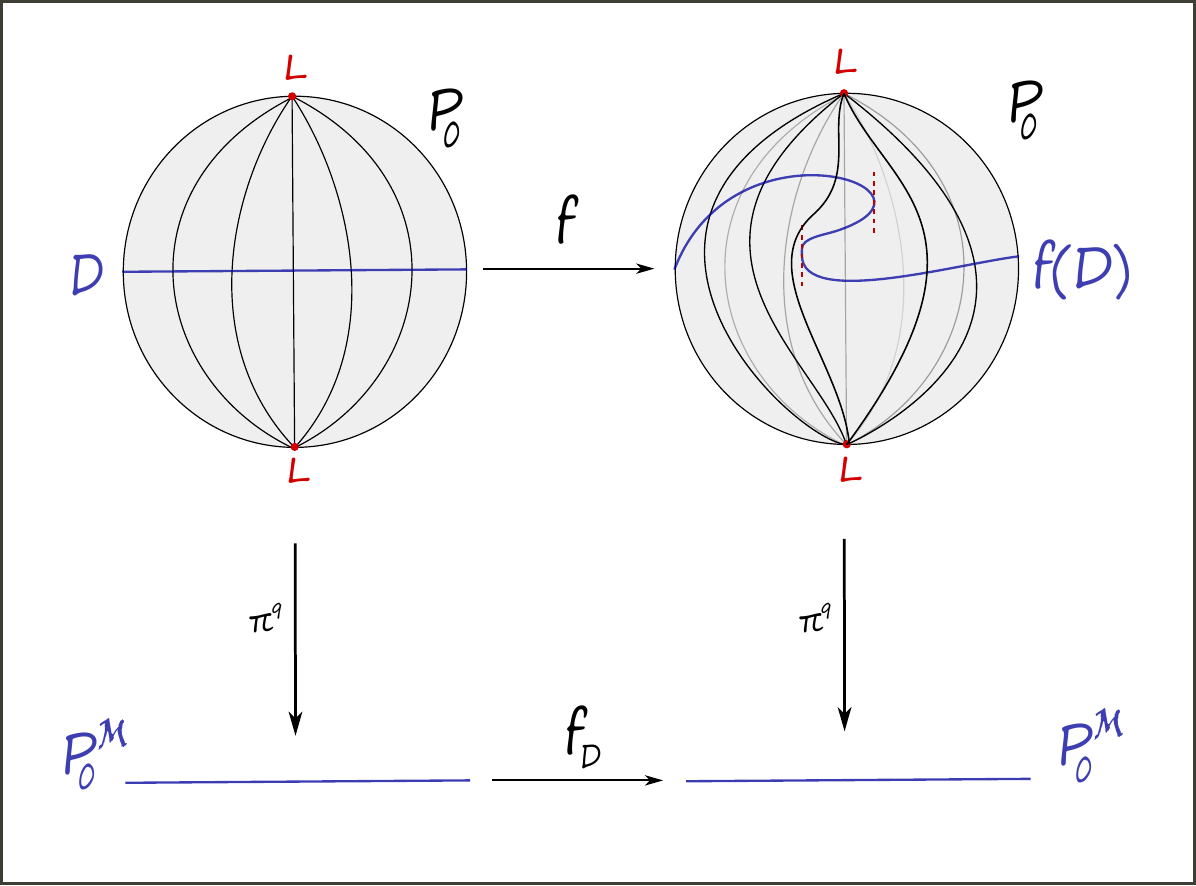}
    \caption{The return map $f_D$ associated to the tomography $D$. Open tangencies of $f(D)$ with the vertical foliation (which could, a priori, theoretically arise from ``foldings'' of the disk $f(D)$) might prevent in general that $f_D$ preserves area. This does \emph{not} happen perturbatively, i.e.\ when we perturb a foliation-preserving map, however. Note that $f(L)=L$.}
    \label{fig:verticaltang}
\end{figure}

For each such symplectic tomography $D$, we have an associated return map 
$$f_D: P_0^\mathcal{M}\rightarrow P_0^\mathcal{M}$$ on the $0$-page of the moduli space, as follows. We identify $P_0^\mathcal{M}$ with $\mathbb{D}_0^2$ via $\mathbb{D}_0^2=(\pi_0 \circ ev \circ \pi_*^{-1})(P_0^\mathcal{M})$, and define $f_D$ by
$$
f_D(u)=\gamma_{s(u)}(\tau(u,D)) \in P_0^\mathcal{M},
$$
where $\tau(u,D)=\min\{t>0: \gamma_{s(u)}(t)\in P_0^\mathcal{M}\}$ is the first return time of the transverse path $\gamma_{s(u)}$ to the $0$-page $P_0^\mathcal{M}$. See Figure \ref{fig:verticaltang}.

The symplectic disk $(D,d\alpha\vert_D)$ is symplectomorphic to $(P_0^\mathcal{M},d\alpha_{\mathcal{M}}\vert_{P_0^\mathcal{M}})$, and both have finite symplectic area. In general, $f_D$ might a priori decrease area. Indeed, $f(D)$ is a symplectic disk in $P_0$ with the same symplectic area as $D$, but it might have an open set of vertical tangencies, i.e.\ intersecting a fiber along an open subset of positive area (as opposed to $D$, which intersects fibers at a single point). Nevertheless, this is not possible for perturbative situations where one perturbs a fiber-wise preserving map, in which case the perturbed $f_D$ still preserves area. 

On the other hand, if $f_D$ is easily seen to be surjective, i.e.\ every holomorphic fiber of $\pi_0$ is intersected by the symplectic disk $f(D)$. Indeed, since $f$ is homotopic to the identity by a smooth homotopy which preserves the boundary, $f(D)$ is homologous to $D$ relative boundary. Moreover, $\partial f(D)$ is a push-off in the Reeb direction of $\partial D$ (they agree if $\partial D$ is an orbit, as observed above), hence they can be homotoped to each other via the Reeb flow of $R_B$, and in particular away from $L$ (the boundary of the holomorphic fibers). It follows that the homological intersection number of $f(D)$ with the fibers agrees with that of $D$, i.e.\ it is $1$. In general, $f_D$ is not necessarily injective. However, this is certainly true in the case where $f$ is close to a fiber-wise preserving map, since otherwise $f(D)$ would have vertical tangencies.  

\medskip

\textbf{Perturbative case.} As observed in the above discussion, if $f$ is sufficiently close to a fiber-wise preserving map, then $f_D$ is an area-preserving homeomorphism of the $2$-disk for every tomography $D$. By Brouwer's translation theorem, we find an interior fixed point for $f_D$; by construction this corresponds to an (interior) fiber-wise $1$-recurrent point in the fixed page $P_0$. Varying vertically the tomography $D$ along $P_0$, we obtain infinitely many such points. If $k\geq 1$, fiber-wise $k$-recurrent points correspond to interior fixed points of the return map $$f_{k,D}(u)=\gamma_{s(u)}(\tau_k(u,D)),$$ where $\tau_k(u,D)$ is the $k$-th return time of the transverse path $\gamma_{s(u)}$ to $P_0^\mathcal{M}$. Note that this map is in general \emph{different} from $f_D^k$; recall that $f_D$ is not the return map of an autonomous flow. Having fixed $k$, using that $D$ and $P_0$ vary in compact families, we can take a sufficiently small perturbation so that $f_{l,D}$ is still an area-preserving homeomorphism for every $l\leq k$, $D$ and any choice of page, and apply the same argument. This finishes the proof of Theorem \ref{thm:application}.

\smallskip

\textbf{Discussion and outlook.} Roughly speaking, one might expect that the results here presented open the possibility of obtaining information on a flow on dimension $5$, from a flow in dimension $3$ (more specifically, on $S^3$). The flow of the holomorphic shadow, by construction, maps holomorphic curves to holomorphic curves, i.e.\ it is a flow in the moduli space. While the ''best approximation'' of the original flow with this property, it potentially alters the dynamics of the original Reeb flow in a significant way. It also forgets dynamical information in the vertical directions, as well as most of the interesting dynamical information at $B$ (it is adapted to study spatial problems rather than planar ones); note that in dimension $3$, the shadow, when seen as a flow on $B$, is just a reparametrization of the original one. Moreover, while the holomorphic shadow admits periodic orbits away from the binding $\mathcal{M}^q_B$ (one or infinitely many, by combining Frank's theorem with Brouwer's translation theorem; and in fact if $\alpha_\mathcal{M}$ is $C^\infty$-generic the union of their images is dense in $S^3$ \cite{Irie}), these may bear no relationship whatsoever with spatial periodic orbits of the original Reeb flow $R_\alpha$; unless the latter vector field satisfies the linearized CR-equation, in which case $R_\alpha=R_\mathcal{M}$, which one expects to happen only for certain integrable systems as the (rotating) Kepler problem. Also, it would be surprising to obtain holomorphic curves which are invariant under the original Reeb flow away from the integrable case (or perhaps perturbations thereof); if such curves existed e.g.\ in the case of the three-body problem where $F=\mathbb{D}^*S^1$ is an annulus and we have twisting at the boundary, one would be able to conclude via the classical Poincar\'e-Birkhoff theorem that there are infinitely many periodic points with arbitrary large period.

On the other hand, one might hope to ``lift'' dynamical properties from the shadow to the original dynamics. For instance, one can envision extracting information related to the \emph{complexity} of the original flow by studying its shadow, at least generically. For instance, from the semi-conjugation (\ref{diag:semiconj}), one obtains that the topological entropy of $\phi_t^{M;\mathcal{M}}$ is no smaller than that of $\phi_t^{S^3;\mathcal{M}}$. Moreover, $C^\infty$-generically, Reeb flows on $S^3$ have positive topological entropy \cite{CDHR}. Since heuristically $\phi_t^{M;\mathcal{M}}$ is 
``less chaotic'' than the original Reeb flow, this suggests that generic adapted Reeb flows on any IP contact $5$-fold might have positive topological entropy, which is moreover generated by spatial orbits (note that if the planar problem has positive entropy, so does the spatial one, since it is an invariant subset). Of course, this is not the case for the integrable limit situations here considered. Note that, the shadow of the $3$-dimensional flow, when viewed as a flow on $B$, has the same topological entropy as the original flow; but, since every homeomorphism of the circle has vanishing entropy, we obtain no interesting information when we semi-conjugate to the moduli space. 

Alternatively, the notion of the transverse shadow, while still forgets dynamical information in the vertical directions (and so in principle it might not be possible to recover the original Reeb flow by knowledge of the collection of all horizontal tomographies for each page), is more reliable than its holomorphic counterpart. The clear disadvantage is that it is not a flow. However, for each choice of symplectic tomography, the failure of the associated transverse shadow to be a (non-autonomous) flow is precisely the existence of fiber-wise intersections. Note that an invariant surface for the return map $f$ is a point $u \in P_0^\mathcal{M}$ which is a fixed point for $f_D$, for \emph{every} symplectic tomography $D$. Note that the notion of a vertical tomography is complicated by the fact that there is no globally defined projection to a fiber.

\appendix

\section{Lorentz metrics and open books}\label{app:cones}

In this appendix, we consider a weaker notion of an open book being adapted to a cone structure as that of Definition \ref{def:stronglyadapted}, and give some curious examples arising from Lorentz metrics.  

\begin{definition}\label{def:weakly}
Consider an everywhere non-trivial cone structure $C$ on a manifold $M$, where $M$ is endowed with an open book $(B,\theta)$. We say that $C$ is weakly adapted to $(B,\theta)$ if there is a Giroux form $\alpha$ for the open book such that: 
\begin{enumerate}
    \item[(1)] The Reeb vector field $R_\alpha$ is positively time-like for $C$; and
    \item[(2)] The contact structure $\xi=\ker \alpha$ is space-like for $C$. 
\end{enumerate}
\end{definition}

Weakly adapted cones arise naturally from Lorentz metrics associated to Giroux forms; see Example (3) below.

\begin{example} We now list some examples.
\begin{enumerate}
    \item  If $\alpha$ is a Giroux form adapted to the open book $(B,\theta)$, then $C=\ker d\alpha$ is a $1$-dimensional cone adapted to $(B,\theta)$, and $\alpha$ is a section of $C$. In fact, every contact form for $\xi=\ker \alpha$ is a section for $C$. The Reeb orbits of $\alpha$ are orbits for the cone, parametrized by $\alpha$.
    \item If $\xi$ is a contact structure supported by the open book $(B,\theta)$, $g$ is a Riemannian metric on a $2n+1$-manifold $M$, $P$ is a page of the open book, and $x\in \mbox{int}(P)$, we can identify the space of directions in the open upper half-space $d_x\theta>0$ with the open upper hemisphere of the $2n$-dimensional unit $g_x$-ball. This space can then be globally parametrized by $2n$ angles $\phi_1,\dots,\phi_{2n}$, and these coordinates are independent on the metric. Given $\varphi=(\varphi_1,\dots,\varphi_{2n}) \in (0,\pi/2)^{2n}$, we let $C_{\varphi}=C\langle K_{\varphi}\rangle$ be the cone structure generated by the collection $K_{\varphi}$ of $1$-dimensional cones $\ker d\alpha$, where $\alpha$ is a contact form for $\xi$ adapted to the open book, and satisfying that $\phi_i(\mathbb{R}^+\cdot\ker d\alpha) \in [\varphi_i,\pi-\varphi_i]$ at every point of every page $P$, for every $i=1,\dots,2n$. Then $C_{\varphi}$ is adapted to $(B,\theta)$. If $\phi_i(\mathbb{R}^+\cdot\ker d\alpha)=\pi$ for $i=1,\dots,2n$, then $C_{\varphi}$ is \emph{centered} at the Reeb vector field $R_\alpha$ of $\alpha$.
    
    \item In the context of general relativity, if $M$ is a time-orientable space-time (i.e.\ it comes endowed with a Lorentz metric $g$, that is a pseudo-metric of signature $(-1,1,\dots,1)$, and a global, continuous, non-vanishing time-like vector field $\partial_t$), then it has associated future and past light cones structures, given by $C_+=\{v \in TM: g(v,v)\leq 0,\; g(v,\partial_t)\leq 0\}, C_-=-C_+=\{v \in TM: g(v,v)\leq 0,\; g(v,\partial_t)\geq 0\}$. Their orbits are the causal curves, i.e.\ respectively the future-directed/past-directed paths.
    
    A compact orientable manifold admits a Lorentz metric if and only if it has vanishing Euler characteristic; in particular, every compact orientable odd-dimensional manifold admits a Lorentz metric. Also, every compact orientable odd-dimensional manifold admits an open book. One can then consider the following situation. 
    Take $(B,\theta)$ an open book for $M$, a Giroux form $\alpha$, and a compatible almost complex structure $J$ on $\xi=\ker \alpha$, chosen so that $\xi\vert_B$ is $J$-invariant. Then $$g_{\alpha,J}=d\alpha(\cdot,J\cdot)-\alpha\otimes \alpha$$ is a Lorentz metric, whose future light cone $C_{\alpha,J}$ is weakly adapted to the open book.    These Lorentz metrics are at the bottom level of the causal hierarchy, as they are totally vicious, i.e.\ they allow for ``time travel''. 
    
    \begin{figure}
    \centering
    \includegraphics[width=1 \linewidth]{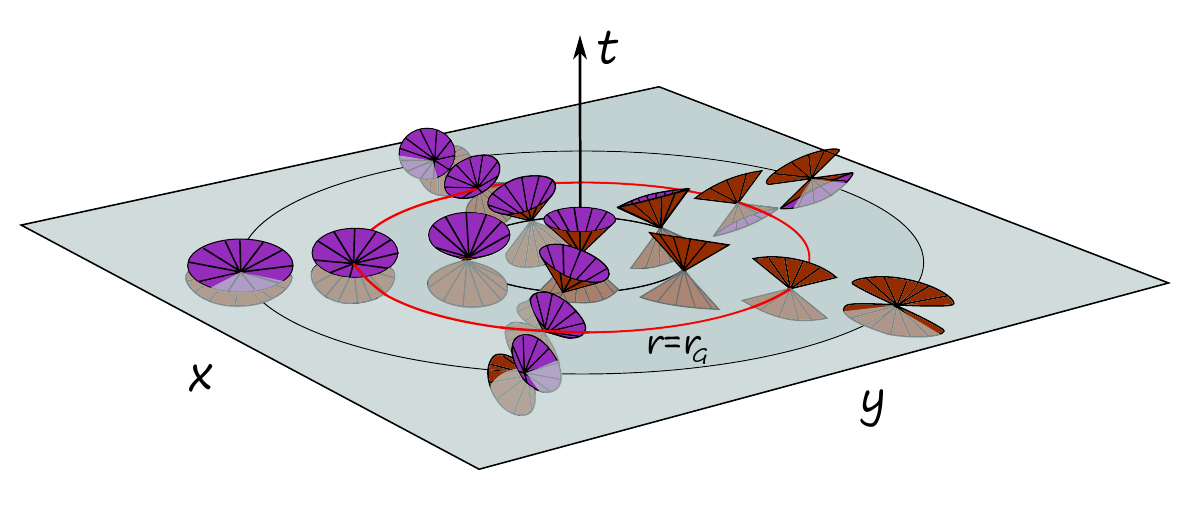}
    \caption{Gödel's metric is static and rotationally symmetric (i.e.\ $\partial_t,\partial_\theta$ are Killing vector fields). Its cone structure tilts as we move radially away from the origin. The circles in the $(x,y)$-plane of constant radius $r$ are space-like for $r<r_G$ less than the value of Gödel's horizon $r_G$, null for $r=r_G$, and time-like for $r>r_G$. The pages of the trivial open book on $\mathbb{R}^3$ in this picture correspond to $t=\mbox{const}$, and the binding, to the circle ``at infinity''. We can clearly guess from this picture that this cone structure should be weakly adapted to the trivial open book on $\mathbb R^3$, where $\alpha=dt-r^2d\theta$ is the standard contact form in these coordinates, with Reeb vector field $\partial_t$. This is indeed true if $r_G$ is chosen adequately. Alternatively, one can turn the picture around by changing charts, and think of the pages as $\theta=\mbox{const}$, the binding as the $t$-axis, $r\leq r_G$ as the collar neighbourhood of the binding, and $r>r_G$ as the mapping torus of the identity monodromy; the contact form $\alpha$ needs to be changed accordingly. Its Reeb vector field is now $\partial_r$ for $r>r_G$, and twists in $r\leq r_G$, coinciding with $\partial_t$ at $r=0$.}
    \label{fig:Godel}
\end{figure}
    
    \item[(4)] As a concrete example of a famous totally vicious Lorentz manifold, satisfying Einstein's field equations, we may consider Gödel's metric (or \emph{Gödel's universe}), given by $(\mathbb{R}^3,g_G)$ with
    $$
    g_G=g_G(a)=-c^2dt^2+\frac{1}{1+(\frac{r}{2a})^2}dr^2+r^2\left(1-\left(\frac{r}{2a}\right)^2\right)d\theta^2+r^2\frac{c}{\sqrt{2}a}dtd\theta,
    $$
    where $c$ is the speed of light, $a>0$ is Gödel's parameter, and $(t,r,\theta)$ are cylindrical coordinates on $\mathbb{R}^3$. We let $r_G=2a$ be the radius of \emph{Gödel's horizon}. Then $\partial_\theta$ is space-like for $r<r_G$, null for $r=r_G$, and time-like for $r>r_G$. Note that $$\lim_{a\rightarrow +\infty}g_G(a)=-c^2dt^2+dr^2+r^2d\theta^2$$ is the metric for flat Minkowski $3$-space. We take the standard contact form on $\mathbb{R}^3$ given by
    $$
    \alpha=dt-r^2d\theta.
    $$
    We clearly have that $R_\alpha=\partial_t$ is time-like, and one readily checks that the restriction of $g_G$ to $\xi=\ker \alpha=\langle \partial_r, \partial_\theta+r^2\partial_t\rangle$ is given by    
    $$
    g_G\vert_\xi=\left(\begin{array}{cc}
        \frac{1}{1+\left( \frac{r}{r_G}\right)^2} & 0 \\
        0 & r^2\left( 1+(\frac{r}{r_G})^2\left(4ac(\sqrt{2}-ac)-1 \right)\right)
    \end{array}\right).
    $$
    Moreover, the expression $4ac(\sqrt{2}-ac)-1$, a degree $2$ polynomial in $a$, is positive if 
    $$
    a \in \left(\frac{\sqrt{2}-1}{2c},\frac{\sqrt{2}+1}{2c}\right).
    $$
    We may e.g.\ take $a=\frac{1}{\sqrt{2}c}$. Given this condition on $a$, $g_G\vert_\xi$ is then a Riemannian metric for every $r>0$, and hence $\xi$ is indeed space-like for $g_G$. See Figure \ref{fig:Godel}. It would be instructive to determine whether $g_G$ on $\mathbb{R}^3$ admits a conformal compactification to a Lorentz metric on $S^3$, which is weakly adapted to the standard open book $S^3=\mathbf{OB}(\mathbb{D}^2,\mathds{1})$; we will not pursue this. 
\end{enumerate}    
\end{example}

\end{document}